\documentclass[a4paper,11pt]{article}
\usepackage[top=3cm,bottom=3cm,left=2.5cm,right=2.5cm]{geometry}
\usepackage[T1]{fontenc}
\usepackage{lmodern,microtype}

\usepackage{xcolor}

\usepackage{subcaption}
\usepackage{tikz,enumerate}
\usetikzlibrary{calc}
\usepackage[numbers,sort&compress]{natbib}

\makeatletter
\renewcommand{\@fnsymbol}[1]{\ifcase#1\or\ensuremath{*}\or\ensuremath{\dagger}\else\@ctrerr\fi}
\makeatother
\usepackage[english]{babel}
\usepackage{amsmath,amsfonts,amsthm,mathtools}
\usepackage{enumitem,booktabs,aliascnt}
\usepackage[hidelinks,hypertexnames=false]{hyperref}
\usepackage[nameinlink,capitalize,noabbrev]{cleveref}
\setlist[enumerate]{leftmargin=2.2em,itemsep=0.15em,topsep=0.25em}
\allowdisplaybreaks
\newcommand{\varnothing}{\emptyset}
\newtheorem{lemma}{Lemma}[section]
\newaliascnt{theorem}{lemma}
\newtheorem{theorem}[theorem]{Theorem}
\aliascntresetthe{theorem}
\newaliascnt{proposition}{lemma}
\newtheorem{proposition}[proposition]{Proposition}
\aliascntresetthe{proposition}
\newaliascnt{corollary}{lemma}
\newtheorem{corollary}[corollary]{Corollary}
\aliascntresetthe{corollary}
\newaliascnt{conjecture}{lemma}
\newtheorem{conjecture}[conjecture]{Conjecture}
\aliascntresetthe{conjecture}
\newaliascnt{problem}{lemma}
\newtheorem{problem}[problem]{Problem}
\aliascntresetthe{problem}

\title{An optimal $\chi$-bound for $\{P_7, \text{diamond}\}$-free graphs}
\author{Kaiyang Lan$^{*1}$ and Hongyang Wang$^{\dagger2}$\\[0.6em]
\small $^{1}$School of Mathematics and Statistics, Minnan Normal University,\\
\small Zhangzhou, Fujian 363000, China\\[0.3em]
\small $^{2}$Department of Mathematics, Hong Kong University of Science and Technology,\\
\small Clear Water Bay, Hong Kong}
\date{}
\begin{document}
\maketitle
\begingroup
\renewcommand{\thefootnote}{\fnsymbol{footnote}}
\footnotetext[1]{Email address: \href{mailto:kylan95@126.com}{kylan95@126.com}.}
\footnotetext[2]{Email address: \href{mailto:hwanghj@connect.ust.hk}{hwanghj@connect.ust.hk} (corresponding author).}
\endgroup
\begin{abstract}
%Schiermeyer~\cite{Schiermeyer2020} asked whether the difference between the chromatic number and the clique number is bounded for graphs with no induced seven-vertex path or diamond. We prove that every such graph $G$ satisfies $\chi(G)\le\omega(G)+3$, where $\chi(G)$ and $\omega(G)$ denote its chromatic number and clique number, respectively.
A \textit{diamond} is a graph obtained by removing an edge from the complete graph on four vertices.
	A graph is $\{P_7, \text{diamond}\}$-\textit{free} if it contains no induced subgraph isomorphic to a seven-vertex path or a diamond.
	Schiermeyer (2020) asked whether there is a constant $c$ such that
	$\chi(G)\leq\omega(G)+c$ for every $\{P_7, \text{diamond}\}$-free graph $G$.
	We answer this question affirmatively by proving \(\chi(G) \leq \omega(G)+3\).
    Moreover, the bound is optimal.
\end{abstract}


\section{Introduction}\label{sec:introduction}

All graphs considered in this paper are finite, simple, and undirected. We write $V(G)$ and $E(G)$ for the vertex set and edge set of a graph $G$. For $X\subseteq V(G)$, let $G[X]$ be the subgraph induced by $X$, let $G-X=G[V(G)\setminus X]$, and write $G-v=G-\{v\}$. For disjoint vertex sets $X,Y$, let $E(X,Y)$ be the set of edges with one end in $X$ and the other in $Y$. The length of a path or cycle is its number of edges.

For a nonnegative integer $k$, a \textit{$k$-coloring} of a graph $G$ is a mapping from $V(G)$ to a set of $k$ colors such that adjacent vertices receive different colors. If such a coloring exists, then $G$ is said to be \textit{$k$-colorable}. The \textit{chromatic number} $\chi(G)$ is the least $k$ for which $G$ is $k$-colorable. A set of pairwise adjacent vertices in a graph \(G\) is called a \textit{clique}, and a set of pairwise nonadjacent vertices in \(G\) is called a \textit{stable set}. The \textit{clique number} \(\omega(G)\) is the maximum size of a clique in \(G\). A graph $G$ is \textit{perfect} if every induced subgraph $H$ of $G$ satisfies $\chi(H)=\omega(H)$. We use $\chi(G)=\omega(G)=0$ for the null graph, whose empty map is a $0$-coloring.

Given a set $\mathcal{H}$ of graphs, we say that $G$ is $\mathcal{H}$-free if $G$ has no induced subgraph isomorphic to any graph in $\mathcal{H}$. If $\mathcal{H}$ consists of a single graph $H$, we abbreviate $\{H\}$-free by $H$-free. A graph class is \textit{hereditary} if it is closed under taking induced subgraphs. A hereditary class $\mathcal{G}$ is \textit{$\chi$-bounded}~\cite{Gyarfas1987} if there is a function $f:\mathbb{N}_0\to\mathbb{N}_0$, where $\mathbb{N}_0=\{0,1,2,\ldots\}$, such that $\chi(H)\le f(\omega(H))$ for every $G\in\mathcal{G}$ and every induced subgraph $H$ of $G$. Such an $f$ is referred to as a \emph{$\chi$-binding function} for $\mathcal{G}$. If $f$ can be chosen to be bounded above by a polynomial in $\omega(H)$, then $\mathcal{G}$ is \textit{polynomially $\chi$-bounded}; equivalently, there is a polynomial $p$ for which $\chi(H)\le p(\omega(H))$ holds for all $G\in\mathcal{G}$ and all induced subgraphs $H$ of $G$. For instance, the class of perfect graphs is polynomially $\chi$-bounded, admitting the identity function $f(x)=x$ as a $\chi$-binding function. While $\chi(G)\ge \omega(G)$ is always true, the reverse inequality is generally false; in fact, triangle-free graphs with arbitrarily large chromatic number are known to exist~\cite{Erdos1959,Mycielski1955,Zykov1949}.

A result of Erd\H{o}s~\cite{Erdos1959} implies that if the class of $H$-free graphs is $\chi$-bounded, then $H$ must be a forest. The Gy\'arf\'as--Sumner conjecture \cite{Gyarfas1973,Sumner1981} proposes the converse statement:

\begin{conjecture}[\cite{Gyarfas1973,Sumner1981}]\label{conjgs}
For every forest $T$, the class of $T$-free graphs is $\chi$-bounded.
\end{conjecture}

We use $P_t$, $C_t$, and $K_t$ to denote, respectively, the path, cycle, and complete graph on $t$ vertices. Gy\'arf\'as \cite{Gyarfas1987} proved Conjecture~\ref{conjgs} for $T=P_t$ with $t\ge2$, showing that every $P_t$-free graph $G$ satisfies $\chi(G)\le (t-1)^{\omega(G)-1}$. Note that this $\chi$-binding function is exponential in $\omega(G)$. Esperet \cite{Esperet2017} conjectured that every $\chi$-bounded class is polynomially $\chi$-bounded. This conjecture was recently disproved by Bria\'nski, Davies, and Walczak \cite{BrianskiDaviesWalczak2024}, who constructed $\chi$-bounded classes that are not polynomially $\chi$-bounded. For $t\ge6$, the question of whether the class of $P_t$-free graphs is polynomially $\chi$-bounded remains open.

The answer is trivial for $t\le 4$, as $P_4$-free graphs are perfect~\cite{Seinsche1974}. For $t=5$, a recent result of Nguyen \cite{Nguyen2025} gives a polynomial bound, improving the earlier exponential bound of Gy\'arf\'as~\cite{Gyarfas1987} and the near-polynomial bound of Scott, Seymour, and Spirkl~\cite{ScottSeymourSpirkl2023}. Huang, Yidong Zhou, and Yuwen Zhou~\cite{HuangZhouZhou2026} recently proved that every nonempty $P_6$-free graph $G$ satisfies the quasi-polynomial bound $\chi(G)\le\omega(G)^{c\log_2\omega(G)}$ for an absolute constant $c$. Apart from these cases, little is known about polynomial $\chi$-bounds for $P_t$-free graphs. Consequently, it is natural to study subclasses of $P_t$-free graphs defined by additional forbidden induced subgraphs. For instance, every $\{P_5,K_3\}$-free graph is 3-colorable, and $C_5$ shows that this bound is tight~\cite{WoegingerSgall2001}; every $\{P_5,C_4\}$-free graph satisfies $\chi(G)\le \lceil \frac{5}{4}\omega(G)\rceil$~\cite{ChoudumKarthickShalu2007}, and this bound is attained by equal-size clique blowups of $C_5$. For an integer $m\ge1$, such a graph is formed by replacing each vertex of $C_5$ with a clique of order $m$, and adding all edges between two replacement cliques if their original vertices are adjacent in $C_5$, and no edges otherwise. For $t=6$, optimal $\chi$-binding functions have been established for the classes of $\{P_6,K_3\}$-free and $\{P_6,C_4\}$-free graphs \cite{BrandstadtKlembtMahfud2006,KarthickMaffray2019}. For $t\ge 7$, Gravier, Ho\`ang, and Maffray \cite{GravierHoangMaffray2003} proved that every $\{P_t,K_3\}$-free graph $G$ satisfies $\chi(G)\le t-2$.

%The currently best known upper bound for $P_7$-free graphs is due to Gravier, Ho\`ang, and Maffray \cite{GravierHoangMaffray2003}, who proved that every $P_7$-free graph $G$ with $\omega(G)\ge 2$ satisfies $\chi(G) \le 5^{\omega(G)-1}$.
%The exponential nature of this bound highlights the difficulty of the general $P_7$-free case and motivates the study of restricted subclasses, such as those with an additional forbidden induced subgraph.
For $P_7$-free graphs, Gravier, Ho\`ang, and Maffray \cite{GravierHoangMaffray2003} proved that every $P_7$-free graph $G$ with $\omega(G)\ge 2$ satisfies $\chi(G) \le 5^{\omega(G)-1}$. More recently, Oum~\cite{Oum2026} lowered the exponential base from five to four, proving $\chi(G)\le3\cdot4^{\omega(G)+4}$. The fact that the best known bound remains exponential highlights the difficulty of the general $P_7$-free case and motivates the study of restricted subclasses, such as those with an additional forbidden induced subgraph. Recently, Huang \cite{Huang2024} showed that every $\{P_7,C_4,C_5\}$-free graph $G$ satisfies $\chi(G)\le \left\lceil\frac{11}{9}\omega(G)\right\rceil$, and this bound is best possible. More recently, Huang, Zhou, and Chang~\cite{HuangZhouChang2026} proved the optimal bound $\chi(G)\le \lceil\frac{5}{4}\omega(G)\rceil$ for every $\{P_7,C_4,C_6\}$-free graph $G$. The $\chi$-boundedness of $P_t$-free graphs and related coloring problems have also been extensively studied; see, for example, \cite{BonomoChudnovskyMaceliSchaudtSteinZhong2018,CameronHuangPenevSivaraman2020,ChenXu2025a,ChudnovskyStacho2018,ChudnovskyMaceliStachoZhong2017,ChudnovskyHuangSpirklZhong2021,ChudnovskyKarthickMaceliMaffray2020,ChudnovskySpirklZhong2024a,ChudnovskySpirklZhong2024b,GaspersHuang2019,HuangKarthick2021,JuJookenGoedgebeurHuang2026,WoegingerSgall2001} and the references therein.

A \textit{diamond} is the graph obtained from $K_4$ by removing an edge. % (see \autoref{fig:threegraphs}).
Randerath~\cite{Randerath1998,SchiermeyerRanderath2019} proved that every $\{P_5, \text{diamond}\}$-free graph \(G\) satisfies $\chi(G) \le \omega(G) + 1$, and this bound is tight. For $\{P_6, \text{diamond}\}$-free graphs, Cameron, Huang and Merkel~\cite{CameronHuangMerkel2021} proved the bound $\chi(G) \le \omega(G) + 3$. Goedgebeur, Huang, Ju, and Merkel~\cite{GoedgebeurHuangJuMerkel2023} subsequently proved the tight bound $\chi(G) \le \max\{\omega(G),6\}$ and gave a polynomial-time algorithm for computing the chromatic number of a graph in this class. For $t\ge4$, Schiermeyer and Randerath~\cite[Corollary~9]{SchiermeyerRanderath2019} proved that every connected $\{P_t, \text{diamond}\}$-free graph $G$ with $\omega(G)\ge2$ satisfies \(\chi(G)\leq (t-2)(\omega(G)-1)\). For \(t=7\), this gives the linear bound \(5\omega(G)-5\), but it leaves open whether the coefficient of \(\omega(G)\) can be reduced all the way to one. Schiermeyer~\cite{Schiermeyer2020} posed the following problem.

\begin{problem}[\cite{Schiermeyer2020}]\label{problem1}
Does there exist a constant $c$ such that every $\{P_7, \text{diamond}\}$-free graph $G$ satisfies $\chi(G)\leq\omega(G)+c$?
\end{problem}

We also note that Problem~\ref{problem1} appears in~\cite{AngeliyaCharKarthick2026}. Two partial results toward Problem~\ref{problem1} are known: Chen, Wu, and Xu~\cite{ChenWuXu2024} proved the bound \(\chi(G) \le \max\{3,\omega(G)\}\) for every \(\{P_7,C_4,\text{diamond}\}\)-free graph \(G\), and Chen and Xu~\cite{ChenXu2025b} established the same bound for every \(\{P_7,C_5,\text{diamond}\}\)-free graph \(G\).

In this paper, we answer Problem~\ref{problem1} affirmatively. Our main result is the following explicit additive bound.

\begin{theorem}\label{thm:intro-main}\label{thm:main}
Every $\{P_7, \text{diamond}\}$-free graph $G$ satisfies \(\chi(G)\leq\omega(G)+3\).
\end{theorem}

\begin{figure}[htbp]
	\centering
	\tikzstyle{v}=[circle, draw, solid, fill=black, inner sep=0pt, minimum width=5pt]
	\begin{tikzpicture}[scale=1.0]
		\node[v] (0) at (0.31537079828623416, 0.11848533500937997) {};
		\node[v] (1) at (-0.8571428571428541, 2.149338558780869) {};
		\node[v] (2) at (-2.029656512571945, 0.11848533500937819) {};
		\node[v] (3) at (-1.7274143833643754, 1.8325856275068177) {};
		\node[v] (4) at (1.7536717215217015, 3.9068840633207014) {};
		\node[v] (5) at (-2.19047619047619, 1.0305390505444965) {};
		\node[v] (6) at (-2.2463282784782983, -3.0213191669548127) {};
		\node[v] (7) at (1.4878844537153215, -0.5584657395811163) {};
		\node[v] (8) at (2.6603981091444133, -1.2354168141716146) {};
		\node[v] (9) at (0.06898075708077211, -1.7490673081032506) {};
		\node[v] (10) at (1.809523809523811, 1.2656416914891166) {};
		\node[v] (11) at (-0.39408105003104277, -0.47681544925168673) {};
		\node[v] (12) at (-2.597685909585894, 2.869734845413758) {};
		\node[v] (13) at (-0.8571428571428559, 3.503240707961862) {};
		\node[v] (14) at (-1.3202046642546685, -0.47681544925168673) {};
		\node[v] (15) at (-4.857142857142858, 1.5007443324337348) {};
		\node[v] (16) at (-0.8571428571428559, 4.857142857142858) {};
		\node[v] (17) at (0.8834001953001831, 2.869734845413758) {};
		\node[v] (18) at (3.1428571428571423, 1.5007443324337348) {};
		\node[v] (19) at (0.47619047619047805, 1.0305390505444942) {};
		\node[v] (20) at (-3.523809523809522, 1.2656416914891166) {};
		\node[v] (21) at (-3.467957435807414, 3.9068840633207014) {};
		\node[v] (22) at (-4.374683823430123, -1.2354168141716146) {};
		\node[v] (23) at (-3.202170168001033, -0.5584657395811199) {};
		\node[v] (24) at (0.013128669078662725, 1.8325856275068162) {};
		\node[v] (25) at (-1.7832664713664839, -1.7490673081032506) {};
		\node[v] (26) at (0.5320425641925874, -3.0213191669548074) {};
		
		\path[draw, blue, thick]
		(0) edge node {} (1) 
		(0) edge node {} (2) 
		(0) edge node {} (3) 
		(0) edge node {} (4) 
		(0) edge node {} (5) 
		(0) edge node {} (6) 
		(0) edge node {} (7) 
		(0) edge node {} (8) 
		(0) edge node {} (9) 
		(0) edge node {} (10) 
		(1) edge node {} (2) 
		(1) edge node {} (11) 
		(1) edge node {} (12) 
		(1) edge node {} (13) 
		(1) edge node {} (14) 
		(1) edge node {} (15) 
		(1) edge node {} (16) 
		(1) edge node {} (17) 
		(1) edge node {} (18) 
		(2) edge node {} (19) 
		(2) edge node {} (20) 
		(2) edge node {} (21) 
		(2) edge node {} (22) 
		(2) edge node {} (23) 
		(2) edge node {} (24) 
		(2) edge node {} (25) 
		(2) edge node {} (26) 
		(3) edge node {} (4) 
		(3) edge node {} (11) 
		(3) edge node {} (12) 
		(3) edge node {} (13) 
		(3) edge node {} (14) 
		(3) edge node {} (19) 
		(3) edge node {} (20) 
		(3) edge node {} (21) 
		(3) edge node {} (22) 
		(4) edge node {} (15) 
		(4) edge node {} (16) 
		(4) edge node {} (17) 
		(4) edge node {} (18) 
		(4) edge node {} (23) 
		(4) edge node {} (24) 
		(4) edge node {} (25) 
		(4) edge node {} (26) 
		(5) edge node {} (6) 
		(5) edge node {} (11) 
		(5) edge node {} (12) 
		(5) edge node {} (15) 
		(5) edge node {} (16) 
		(5) edge node {} (19) 
		(5) edge node {} (20) 
		(5) edge node {} (23) 
		(5) edge node {} (24) 
		(6) edge node {} (13) 
		(6) edge node {} (14) 
		(6) edge node {} (17) 
		(6) edge node {} (18) 
		(6) edge node {} (21) 
		(6) edge node {} (22) 
		(6) edge node {} (25) 
		(6) edge node {} (26) 
		(7) edge node {} (8) 
		(7) edge node {} (11) 
		(7) edge node {} (13) 
		(7) edge node {} (15) 
		(7) edge node {} (17) 
		(7) edge node {} (19) 
		(7) edge node {} (21) 
		(7) edge node {} (23) 
		(7) edge node {} (25) 
		(8) edge node {} (12) 
		(8) edge node {} (14) 
		(8) edge node {} (16) 
		(8) edge node {} (18) 
		(8) edge node {} (20) 
		(8) edge node {} (22) 
		(8) edge node {} (24) 
		(8) edge node {} (26) 
		(9) edge node {} (10) 
		(9) edge node {} (11) 
		(9) edge node {} (14) 
		(9) edge node {} (16) 
		(9) edge node {} (17) 
		(9) edge node {} (20) 
		(9) edge node {} (21) 
		(9) edge node {} (23) 
		(9) edge node {} (26) 
		(10) edge node {} (12) 
		(10) edge node {} (13) 
		(10) edge node {} (15) 
		(10) edge node {} (18) 
		(10) edge node {} (19) 
		(10) edge node {} (22) 
		(10) edge node {} (24) 
		(10) edge node {} (25) 
		(11) edge node {} (18) 
		(11) edge node {} (22) 
		(11) edge node {} (24) 
		(11) edge node {} (25) 
		(11) edge node {} (26) 
		(12) edge node {} (17) 
		(12) edge node {} (21) 
		(12) edge node {} (23) 
		(12) edge node {} (25) 
		(12) edge node {} (26) 
		(13) edge node {} (16) 
		(13) edge node {} (20) 
		(13) edge node {} (23) 
		(13) edge node {} (24) 
		(13) edge node {} (26) 
		(14) edge node {} (15) 
		(14) edge node {} (19) 
		(14) edge node {} (23) 
		(14) edge node {} (24) 
		(14) edge node {} (25) 
		(15) edge node {} (20) 
		(15) edge node {} (21) 
		(15) edge node {} (22) 
		(15) edge node {} (26) 
		(16) edge node {} (19) 
		(16) edge node {} (21) 
		(16) edge node {} (22) 
		(16) edge node {} (25) 
		(17) edge node {} (19) 
		(17) edge node {} (20) 
		(17) edge node {} (22) 
		(17) edge node {} (24) 
		(18) edge node {} (19) 
		(18) edge node {} (20) 
		(18) edge node {} (21) 
		(18) edge node {} (23) 
		(19) edge node {} (26) 
		(20) edge node {} (25) 
		(21) edge node {} (24) 
		(22) edge node {} (23);
	\end{tikzpicture}
	\caption{The complement of the Schl\"afli graph from House of Graphs \#19273.}	
	%	\caption{The 27-vertex 9-regular graph from House of Graphs \#19273.}
	\label{fig:house19273}
\end{figure}
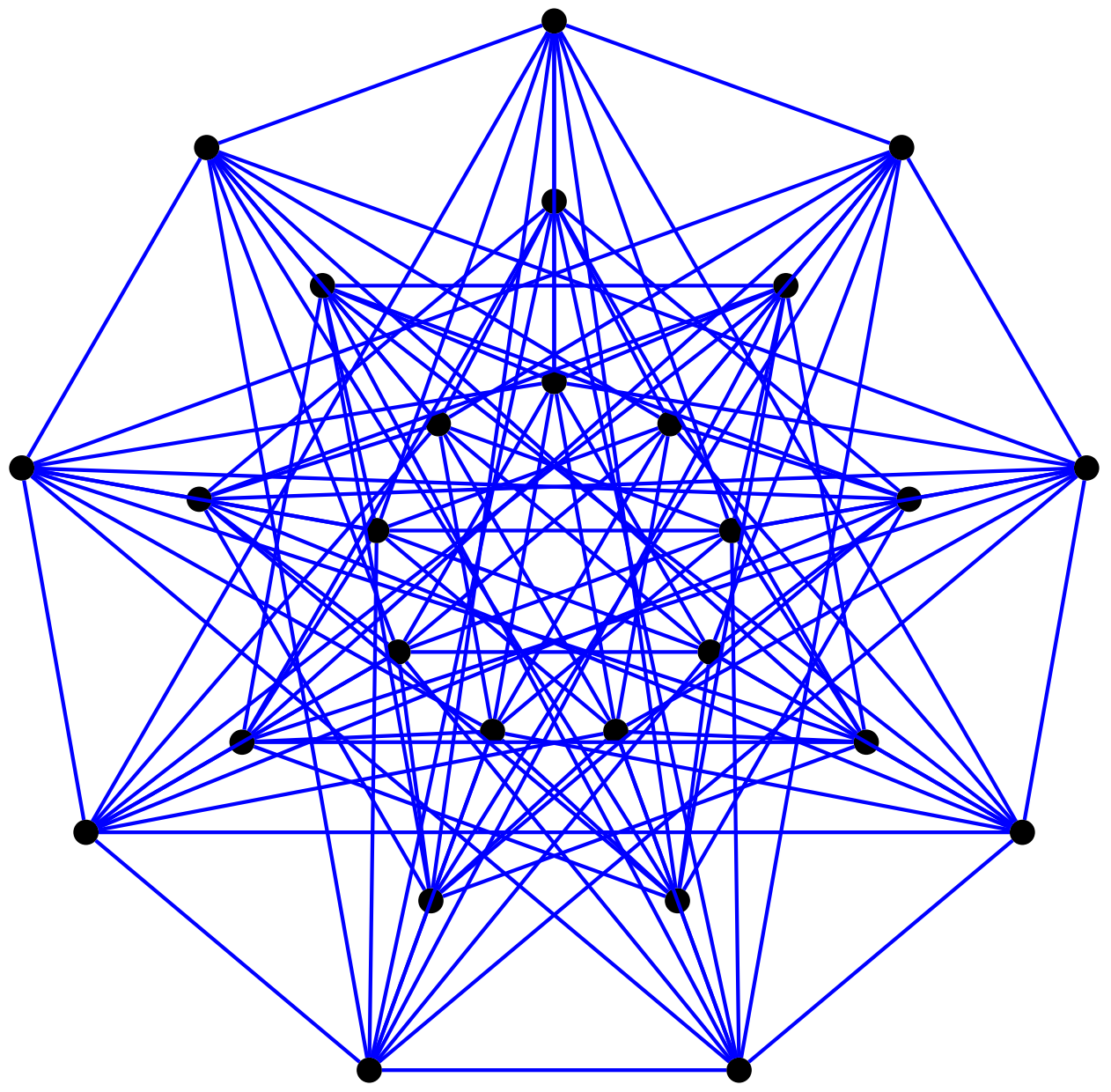

The additive constant three in Theorem~\ref{thm:intro-main} is best possible, as shown by the complement of the $27$-vertex Schl\"afli graph~\cite{CameronHuangMerkel2021}\footnote{\url{https://houseofgraphs.org/graphs/19273}} (see \autoref{fig:house19273}). Moreover, Goedgebeur, Huang, Ju and Merkel~\cite{GoedgebeurHuangJuMerkel2023} gave a $25$-vertex $\{P_6, \text{diamond}\}$-free graph with chromatic number six and clique number three. For every integer $n\ge25$, adding $n-25$ isolated vertices to this graph gives an $n$-vertex $\{P_7, \text{diamond}\}$-free graph $G$ with $\chi(G)=6=\omega(G)+3$.
%Theorem~\ref{thm:intro-main} naturally leads to the following problem.
A natural next step is to ask whether an additive-constant bound still holds when $P_7$ is replaced by a longer path. This motivates the following problem.

\begin{problem}
	For every integer $t\ge 8$, does there exist a constant $c_t$ such that every $\{P_t, \text{diamond}\}$-free graph $G$ satisfies $\chi(G) \le \omega(G) + c_t$?
\end{problem}

%Let \(c_7\) be the least integer \(c\), if it exists, for which every $\{P_7, \text{diamond}\}$-free graph \(G\) satisfies \(\chi\leq\omega(G)+c\).
%The complement of the Schl\"afli graph\footnote{\url{https://houseofgraphs.org/graphs/19273}} (see \autoref{fig:house19273}) is \(P_6\)-free and hence \(P_7\)-free, so the result of Cameron, Huang and Merkel~\cite{CameronHuangMerkel2021} gives \(c_7\geq3\).
%Theorem~\ref{thm:intro-main} therefore yields \(3\leq c_7\leq13\).
%We do not expect the constant \(13\) to be optimal; rather, we would expect $c_7 = 3$.

\section{Preliminaries and rooted colorings}\label{sec:preliminaries}

Let $G$ be a graph. For a vertex \(v\in V(G)\), the neighborhood of $v$ is $N_G(v)$, its degree is $d_G(v)=|N_G(v)|$, and $\delta(G)$ is the minimum degree. We omit the graph subscript when it is clear. For a set \(X\subseteq V(G)\), we write $N_X(v)=N(v)\cap X$, $N(X)=\bigcup_{x\in X}N(x)$ and $N[X]=X\cup N(X)$. The \textit{external boundary}, or \textit{full boundary}, of $X$ within an induced graph $J$ is $N(X)\cap(V(J)\setminus X)$. In particular, $N(X)$ itself may meet $X$. Relative to an induced graph $J$, the zero graph of $X$ is $G[V(J)\setminus N[X]]$; its components are the zero components of $X$. An induced subgraph of this zero graph is a zero subgraph of $X$. We write $xy=1$ or $xy=0$ for adjacency or nonadjacency. Two disjoint vertex sets of $G$ are \textit{complete} to each other if every vertex of one is adjacent to every vertex of the other, and \textit{anticomplete} if there is no edge between them. A vertex \textit{sees} another vertex if they are adjacent, \textit{completes} a set if it is adjacent to every vertex in it, and \textit{misses} a set if it has no neighbor in it. Vertex sets in coloring statements denote their induced subgraphs; for $v\notin V(J)$, write $J+v=G[V(J)\cup\{v\}]$. A \textit{cluster graph} is a disjoint union of cliques, called its blocks. A matching is a set of pairwise vertex-disjoint edges, and it is perfect if it covers every vertex under consideration; the ends of a matching edge are mates. The phrase ``a matching and isolates'' means a graph of maximum degree at most one. A vertex cover is a set meeting every edge. A stable transversal of a clique family is a stable set meeting every member of the family. A set dominates another disjoint set when every vertex of the latter has a neighbor in the former. A clique cutset is a clique whose deletion leaves at least two nonempty components. The null graph has $\chi=\omega=0$; a nonnull edgeless graph has $\chi=\omega=1$. A vertex outside a clique has an \textit{owner} there when it has exactly one neighbor, and is \textit{unowned} when it has none. Given a coloring of \(G\), a \textit{palette} is a set of available colors; a \textit{color class} consists of all vertices receiving a specified color. Palette colors may be unused, in which case their color classes are empty.

\subsection{Ownership and criticality}

\begin{lemma}
\label{lem:clique-ownership}

If $Q$ is a maximal clique of size at least two in a diamond-free graph, every vertex outside $Q$ has at most one neighbor in $Q$. Distinct maximal cliques intersect in at most one vertex and, when they intersect, their exclusive parts are anticomplete. Moreover, every neighborhood induces a disjoint union of cliques.

\end{lemma}
\begin{proof}

A vertex outside $Q$ cannot be complete to $Q$, by maximality. Two neighbors and a nonneighbor in $Q$, together with that vertex, would induce a diamond. The intersection assertions follow because a vertex outside one clique already sees any shared vertex. Finally, an induced three-vertex path in a neighborhood, together with the vertex defining the neighborhood, would induce a diamond. A graph without an induced three-vertex path is a disjoint union of cliques.
\end{proof}

For $k\ge1$, a graph $G$ is \emph{$k$-vertex-critical} if $\chi(G)=k$ and $\chi(G-v)<k$ for every $v\in V(G)$. Equivalently, $\chi(G)=k$ and every proper induced subgraph is $(k-1)$-colorable. Unless specified otherwise, critical means vertex-critical.

\begin{lemma}
\label{lem:critical-core}

For every integer $1\le k\le\chi(G)$, the graph $G$ has an induced $k$-vertex-critical subgraph. If $H$ is $k$-vertex-critical, then $H$ is connected and $\delta(H)\ge k-1$, and $H$ has no clique cutset. For every nonempty stable set $I$, $\chi(H-I)=k-1$. No distinct vertices $x,y$ satisfy $N_H(x)\subseteq N_H(y)$.

\end{lemma}
\begin{proof}

Choose an inclusion-minimal $X\subseteq V(G)$ with $\chi(G[X])\ge k$. Deleting one vertex lowers chromatic number by at most one, so $\chi(G[X])=k$ and $G[X]$ is vertex-critical. Criticality gives connectedness by taking a maximum-chromatic component, and $\delta(H)\ge k-1$ by extending a $(k-1)$-coloring of $H-v$ when $d_H(v)\le k-2$. For a clique cutset, $(k-1)$-color the proper induced graphs formed by the cutset and each component of its deletion. Permute their colors to agree on the cutset and combine the colorings, contradicting $\chi(H)=k$. Criticality and one new color on $I$ give $\chi(H-I)=k-1$. If $N_H(x)\subseteq N_H(y)$, then $xy=0$, and a $(k-1)$-coloring of $H-x$ extends by giving $x$ the color of $y$.
\end{proof}

A \emph{critical core} of $G$ is an induced vertex-critical subgraph with chromatic number $\chi(G)$; its clique number can be smaller than $\omega(G)$.

\subsection{Rooted colorings}
\label{sec:rooted}
A vertex contacts $X$ if it has a neighbor in $X$. An exterior vertex $v\notin X$ is a root on $X$, with profile $N_X(v)$. A rooted $P_t$ into $X$ is an induced path $v,x_1,\ldots,x_{t-1}$ with all $x_i\in X$. The root is good on $X$ if it starts no such $P_4$, and bad otherwise. Missing $X$ implies goodness, but the coloring conclusions require contact.

Two path observations will be used. If an exterior vertex $v$ contacts connected induced $J$ and starts no $P_5$ into it, each component of $J-N_J(v)$ has nonempty full $J$-boundary in $N_J(v)$, consisting of good roots: prefix $v$ to a rooted $P_4$. No vertex in a $P_7$-free graph is bad on two anticomplete sets: its two rooted $P_4$'s concatenate to a $P_7$.

\begin{lemma}
\label{lem:q4-good-extension}
\label{lem:good}

Let the ambient graph be diamond-free with clique number at most $q\ge3$. Let $D$ be connected and induced, contacted by an exterior vertex $v$ good on $D$, and put $B=N_D(v)$. Every component of $D-B$ is a clique of order at most $q-1$. A vertex of $B$ contacting such a component completes it, and a nonsingleton component has exactly one neighbor in $B$. The graph $D+v$ has a $q$-coloring in which $B$ uses only the first $q-1$ colors.

For $q=3$, two further conclusions hold. If $D$ is triangle-free, then $B$ is a vertex cover inducing a matching and isolated vertices (a \emph{matching vertex cover}), and $D-B$ is stable. If $T$ is a triangle in $D$ and $D\ne T$, exactly one vertex $t\in T$ has a neighbor in $D-T$, and every exterior good root contacting $D$ contacts $T$ at $t$.

\end{lemma}
\begin{proof}

If $b\in B$ contacts a component $Z$ of $D-B$, it completes $Z$; otherwise an adjacency transition gives an induced $v,b,z,z'$. Connectedness supplies such a $b$, whose cluster neighborhood forces $Z$ to be a clique. If $|Z|\ge2$, its attachment is unique: two nonadjacent attachments give a diamond on an edge of $Z$, while two adjacent attachments give a diamond with $v$ and a vertex of $Z$. Thus $|Z|\le q-1$.

The cluster graph $B$ has blocks of order at most $q-1$. Color them with $q-1$ colors, and give $v$ and all singleton components of $D-B$ one additional color. Color each nonsingleton $Z$ with the $q-1$ colors different from its unique attachment's color. Distinct components of $D-B$ are anticomplete, so the coloring is proper.

For $q=3$, triangle-freeness excludes an edge component with its attachment. For a triangle $T\subseteq D$, any good root contacting $D$ must contact $T$, since its missed components have order at most two, and has at most one neighbor on $T$ by \cref{lem:clique-ownership}. The two missed triangle vertices form an edge component with unique attachment, so have no $D$-neighbors outside $T$. Connectedness gives the asserted common contact.
\end{proof}

In a diamond-free graph of clique number at most three, every neighborhood induces a matching and isolated vertices, a vertex outside a triangle has at most one neighbor on it, and every edge has at most one common neighbor, by \cref{lem:clique-ownership}. We call these consequences triangle sparsity.

\subsection{Crossed contacts and shared attachments}

The interface of two disjoint vertex sets is the entire edge set between them.

\begin{lemma}
\label{local-crossed-contacts}
\label{lem:crossed-contact}
Let $B,C$ be disjoint cliques in a $P_7$-free graph, let $r\notin B\cup C$ miss $B\cup C$, and let $p,q\in N(r)\setminus(B\cup C)$ each have at most one neighbor in each clique. Suppose either

\begin{enumerate}[label=(\arabic*)]
\item $B,C$ are anticomplete and both have order at least three; or
\item their interface is a matching of size at most two, both have order
at least three, $q$ misses $B$ if $|B|=3$, and $p$ misses $C$ if $|C|=3$.

\end{enumerate}
If $pb=qc=1$ and $pc=qb=bc=0$ for $b\in B,c\in C$, then $pq=1$.

\end{lemma}
\begin{proof}
Suppose $pq=0$ and put
\[
X=B\setminus(\{b\}\cup N_B(q)\cup N_B(c)),\quad
Y=C\setminus(\{c\}\cup N_C(p)\cup N_C(b)).
\]
Both sets are nonempty. In (1) at most two vertices are excluded from each clique. In (2) at most three are excluded, and a clique of order three loses at most two by its extra missed-profile condition.

In (1) every $X$--$Y$ pair is a nonedge. In (2), if every such pair were an edge, the matching would force $X=\{b'\},Y=\{c'\}$. The same size bounds show that $N_B(c)$ and $N_C(b)$ are nonempty: without either term the corresponding candidate set has at least two vertices. Those two interface edges and $b'c'$ are distinct, by $bc=0$ and the candidate exclusions, contradicting the size bound. Thus choose $b'\in X,c'\in Y$ with $b'c'=0$.

Now $b',b,p,r,q,c,c'$ is induced. The profiles remove unwanted contacts with $p,q$, the center misses both cliques, and the four selected $B$--$C$ pairs are nonedges. The crossed profiles give $p\ne q$; the remaining distinctness follows from disjointness, adjacency to $r$, and the candidate definitions.
\end{proof}

\begin{lemma}
\label{lem:q4-shared-boundaries}

In a $\{P_7,\mathrm{diamond}\}$-free graph, let $A$ induce a cluster graph. Let $r$ complete $A$ and miss two anticomplete cliques $E,F$ of order at least three, and suppose each vertex of $A$ has at most one neighbor in each clique. If every vertex of $E\cup F$ has an $A$-neighbor, some block $B$ of $A$ dominates $E\cup F$. In particular, this is impossible if $A$ is a matching and isolates. Under this coverage hypothesis, if every block of $A$ has order at most three, then $B,E,F$ are triangles, $B$ is perfectly matched to each of $E,F$, and $A\setminus B$ misses both.

\end{lemma}
\begin{proof}

If $p\in A$ sees $e\in E$ and misses $f\in F$, while $q\in A$ sees $f$ and misses $e$, then $pq$ is an edge by \cref{local-crossed-contacts}, using its anticomplete case. This rule requires no coverage hypothesis.

Assume coverage and put $E_B=N_E(B)$, $F_B=N_F(B)$ for each block $B$. If $e\in E\setminus E_B$ and $q\in A$ sees $e$, the rule forces $F_B\subseteq N_F(q)$, so $|F_B|\le1$. By symmetry and $|E|,|F|\ge3$, a block meeting two positions of either clique dominates both, and has order at least $\max\{|E|,|F|\}$, excluding matching $A$.

Suppose no such block exists. Represent each block contacting $E\cup F$ by its single $E$- and $F$-contacts, using $\bot$ for absence. Two block profiles differing in both coordinates contain a real crossed pair, contradicting the rule. Thus two different $E$ coordinates share one $F$ coordinate; any other $F$ coordinate must then share both $E$ coordinates, impossible under coverage.

When blocks have order at most three, $B,E,F$ are triangles, with bijective contacts. A root in $A\setminus B$ seeing $e\in E$ misses two $F$ positions. The rule forces their distinct $B$-neighbors both to see $e$, impossible. Symmetrically no root in $A\setminus B$ sees $F$.
\end{proof}

\subsection{A minimal counterexample}
\label{sec:universal-reduction}
\begin{lemma}
\label{lem:minimal-counterexample}

If the bound $\chi\le\omega+3$ fails in the $\{P_7,\mathrm{diamond}\}$-free class, there is a counterexample $G$ with $q=\omega(G)\le4$ that is $(q+4)$-vertex-critical, has minimum degree at least $q+3$, and has no stable set meeting all maximum cliques. If $q=4$, two disjoint maximum cliques have at least three edges between them.

\end{lemma}
\begin{proof}

Minimize the clique number of a counterexample, then its order. Every proper induced subgraph $H$ satisfies $\chi(H)\le\omega(H)+3$. Deleting and restoring a vertex gives $\chi(G)=q+4$ and vertex-criticality; \cref{lem:critical-core} gives the degree bound. A stable set meeting all maximum cliques would leave a $(q+2)$-colorable graph, contradicting $\chi(G)=q+4$ after one new color is added.

Suppose $q\ge5$, or $q=4$ and every two disjoint maximum cliques have at most two cross edges. Choose a stable set $S$ meeting as many maximum cliques as possible, and let $Q=\{q_1,\ldots,q_q\}$ be missed. The sets $S_i=S\cap N(q_i)$ are disjoint and nonempty. Extend $(S\setminus S_i)\cup\{q_i\}$ to a maximal stable set $T_i$. It misses a maximum clique $R_i$ met by $S$, whose sole $S$-vertex $s_i$ belongs to $S_i$. Moreover $R_i\cap Q=\varnothing$: a common vertex could only be $q_i$, which $T_i$ contains.

Fix $j$ and write $R=R_j$, $s=s_j$, $T=T_j$. For $t\in R\setminus\{s\}$, restricted to vertices missing $Q$ when $q=4$, choose $z\in T$ adjacent to $t$. Ownership makes $s,t,z$ induced and gives $tq_j=0$. Since $q_j\in T$, also $z\notin Q$. Let $o,k$ be the possible $Q$-owner indices of $t,z$, omitting nonexistent indices. We claim $ts_i=1$ whenever $i\notin\{j,o,k\}$.

Otherwise put $A=R_i$ and $h=s_i$. Here $h\in S\cap T$ has $Q$-owner $i$, so $h\ne s,z$. The assumed nonedge and stability make $h$ miss $s,t,z$, so the path lies outside $A$. The owner indices, ownership in $R$, and stability of $T$ make $q_i,q_j$ miss $t,z$. Every vertex of $A-h$ contacts $\{q_j,s,t,z\}$, or it prefixes $h,q_i,q_j,s,t,z$ to give $P_7$. If $q\ge5$, ownership forces $q=5$ and four distinct contacts $c_j,c_s,c_t,c_z$ in $A-h$. Choose $q_\ell$ with $\ell\notin\{i,j,o,k\}$. It misses $h,s,t,z$ and at least one of $c_s,c_z$. Accordingly $q_\ell,q_i,h,c_s,s,t,z$ or $q_\ell,q_i,h,c_z,z,t,s$ is induced, a contradiction.

For $q=4$, each path vertex has an $A$-owner, with consecutive owners distinct. Otherwise it covers exactly two positions of $A-h$, and $q_j$ contacts the third. An endpoint has an owner $c$ missed by the other two path vertices. Choose $q_\ell\in Q\setminus\{q_i,q_j\}$ missing $z$. The two $Q$--$A$ contacts are exhausted, so $q_\ell,q_i,h,c$ followed by the path from that endpoint is induced. Thus the consecutive owners are distinct, which gives $A\cap R=\varnothing$. The contacts of $s,t$ exhaust the $A$--$R$ interface. Choose $p\in R\setminus\{s,t\}$ missing $Q$; this is possible since $sq_j$ uses one of its at most two contacts. Then $p$ misses $A,Q,z$. Let $c$ own $z$ in $A$ and choose $x\in Q\setminus\{q_i\}$ missing $z,c$. Now $p,t,z,c,h,q_i,x$ is induced. Ownership and the established nonedges exclude every chord in these paths. This proves the claim.

Stability puts every $s_i$ with $i\ne j$ outside $R$, with at most one neighbor there. For $q\ge5$, the claim gives each of its $q-1$ vertices other than $s$ at least $q-3\ge2$ such neighbors, impossible. For $q=4$, at least two vertices of $R-s$ miss $Q$, because $sq_j$ uses one of at most two cross edges. Each sees at least two of the three retained $s_i$, again impossible.
\end{proof}

\begin{lemma}
\label{lem:triangle-free-five}

Every triangle-free $P_7$-free graph is five-colorable.

\end{lemma}
\begin{proof}

Suppose a connected triangle-free graph $G$ has chromatic number at least six. Let $C_0=G$ and choose $v_0\in V(G)$. For $i=0,\ldots,4$, maintain $\chi(C_i)\ge6-i$. Deleting the stable set $N_{C_i}(v_i)$ leaves a component $C_{i+1}$ with chromatic number at least $5-i$. Connectedness of $C_i$ supplies $v_{i+1}\in N_{C_i}(v_i)$ contacting $C_{i+1}$. Since $\chi(C_1)\ge5$, the first component excludes $v_0$, which is isolated after deletion. Inductively, $C_{i+1}$ is anticomplete to $v_0,\ldots,v_i$, so the chosen vertices form an induced path. Append a neighbor $v_6$ of $v_5$ in $C_5$ to obtain an induced $P_7$. Reuse five colors on distinct components to finish.
\end{proof}

\subsection{A root with a matching neighborhood}

\begin{lemma}
\label{lem:q4-matching-root}
\label{thm:rooted-four}
\label{lem:nonbipartite-boundary-cores}
\label{lem:q4-zero-four-clique}

Let $C$ be connected and induced in a $\{P_7,\mathrm{diamond}\}$-free graph of clique number at most four. Let an exterior root $r$ contact $C$ and start no $P_5$ into $C$. Every induced $J\subseteq C$ for which $N_J(r)$ is a matching and isolates is four-colorable. The neighborhood condition is automatic when the ambient clique number is at most three.

\end{lemma}
\begin{proof}
The last assertion follows from \cref{lem:clique-ownership}. If some such $J$ is not four-colorable, \cref{lem:critical-core} supplies an induced five-vertex-critical $H\subseteq J$. The root contacts $H$: otherwise $H$ is in a component $Z$ of $C-N_C(r)$, whose good attachment supplied by \cref{sec:rooted} four-colors $Z$ by \cref{lem:q4-good-extension}. Put $A=N_H(r)$. Call the components of $H-A$ zero components.

Criticality gives $\delta(H)\ge4$. Every zero component has a nonempty good boundary in $A$: prefix $r$ to any boundary-root $P_4$. Adjacent $A$-vertices have disjoint zero profiles, by the diamond with $r$; hence every attachment group $N_A(f)$ is stable. On a clique of order at least three an $A$-vertex has zero, one, or all contacts, by diamond-freeness. Call a single-contact root \emph{partial} and a completing root a \emph{completer}. A triangle has at most one completer (two give a diamond or a five-clique), and a four-clique has none. Each position of an exact zero triangle has at least two $A$-neighbors, so retains a partial root after deletion of its completer; every zero four-clique position has a partial root.

\emph{Boundary colors.} For every nonbipartite triangle-free zero component $D$, put $S_D=N_A(D)$, its nonempty full $H$-boundary. For each $p\in S_D$, the graph $D+p$ has clique number at most three; \cref{lem:q4-good-extension} makes its profile a matching vertex cover with stable complement.

The boundary $S_D$ is stable: adjacent roots have disjoint profiles by the diamond with $r$, so their two vertex covers would bipartition $D$.

For boundaries $S=N_A(D)$ and $T=N_A(E)$ of distinct such components, suppose $ab$ is an edge with $a\in S$ and $b\in T$. Stability gives $a\notin T$ and $b\notin S$. If $c\in S\setminus(T\cup\{a\})$, choose transitions $c,x,y$ in $D$ and $b,z,w$ in $E$, with $cy=bw=0$; the matching profiles are proper subsets of their connected nonbipartite components. Then $y,x,c,r,b,z,w$ is an induced $P_7$, by full boundaries, component separation, and the matching on $A$. Hence $S-\{a\}\subseteq T$. Reflection gives $T-\{b\}\subseteq S$, and therefore $T=S-\{a\}+\{b\}$.

Let $A_*$ be the union of these boundaries, and let $W_0$ consist of the isolated vertices of $H[A_*]$. Then $W_0$ is stable. Suppose a boundary $S$ loses distinct vertices $a,c$ from $W_0$, with respective mates $b,d\in A_*$. Choose boundaries $T,U$ containing $b,d$. The exchange rule gives $T=S-\{a\}+\{b\}$ and $U=S-\{c\}+\{d\}$. Since $b\in T$, $a\in U$, and $ab$ is an edge, applying it again gives $U=T-\{b\}+\{a\}=S$, contradicting $d\notin S$. Thus $|S_D\setminus W_0|\le1$ for every such $D$.

We use the following rules to color $H$. The \emph{crossing rule} of \cref{lem:q4-shared-boundaries} applies to anticomplete cliques $E,F$ of order at least three and roots in $B\subseteq A$ having at most one contact in each: $pe=qf=1$, $pf=qe=0$ force $pq=1$; moreover, $B$ cannot dominate both cliques. The \emph{palette rule} starts with a proper $1,2$-coloring of $A$. Bipartite zero components use $3,4$. A nonbipartite triangle-free component $D$ also uses $2,3,4$ if at most one boundary root has color $2$: choose that root as $p$, or any root if none does; color the matching $N_D(p)$ with $3,4$ and its stable complement with $2$, using the boundary-cover property above. This respects every boundary color. A zero triangle uses $2,3,4$ if some position has its entire attachment group colored $1$: give that position color $2$.

\emph{Zero shapes.} Suppose a nonclique zero component $D_0$ contains a maximal clique $T$ of order at least three, maximal within $D_0$. A good boundary root must contact $T$, or the missed clique containing $T$ has a complete attachment enlarging $T$, by \cref{lem:q4-good-extension}. If its sole contact is $t$, the missed component containing $T-t$ has unique attachment $t$; maximality makes that component exactly $T-t$. Thus $T-t$ misses $D_0-T$. Different singleton contact positions would isolate $T$ from $D_0-T$, so all partial roots use the same $t$. Vertices of $T-t$ then have degree at most three in $H$: only their $T$-neighbors and, for $|T|=3$, the possible completer. This is impossible. If no partial root occurs, the nonempty boundary is a single completer. A four-coloring of $H-D_0$ extends across $D_0$ by \cref{lem:q4-good-extension}, agreeing at that root. Hence every triangle-containing zero component is a triangle or a four-clique.

\emph{Mixed components.} Let $D$ be a nonbipartite triangle-free zero component. Its full boundary $S=N_A(D)$ is stable with profiles that are matching vertex covers, as shown above. We have $|S|\ge3$: a four-coloring of $H-D$ extends if the boundary uses one color, by reserving it and three-coloring $D$ using \cref{lem:q4-good-extension}. With two distinct boundary colors $3,4$, choose the color-$3$ root $p$ and color $N_D(p)$ with $1,2$, its stable complement with $3$.

Suppose a zero clique $E$ of order at least three coexists with $D$. The partial $E$-contacts of $S$ occupy at most two positions. Otherwise choose $p,q,s\in S$ with distinct singleton contacts $f_p,f_q,f_s\in E$. A shortest odd cycle $Z\subseteq D$ is induced and has length five or seven: triangle-freeness excludes three, and seven consecutive vertices of a longer induced odd cycle give $P_7$. Label its vertices by the roots in $\{p,q,s\}$ seeing them. Consecutive labels have union $\{p,q,s\}$ and intersection of size at most one, by edge coverage and the diamond on a cycle edge with two roots. No root occurs in three consecutive labels, since its profile is a matching.

Some singleton label is contained in a doubleton label. Otherwise, if $\{s\}$ occurs, every doubleton is $\{p,q\}$; neither two labels of size at most one nor two of size at least two can be consecutive, by the union and intersection conditions. Their sizes would bipartition $Z$. With no singleton, an empty label forces alternating empty and full labels. If neither occurs, all labels are doubletons; consecutive ones are distinct, and their common root must be absent from the next, forcing the third doubleton. Period three cannot close on a cycle of length five or seven.

Rename so $x,y$ have labels $\{q,s\},\{s\}$. Both miss $p$, so $xy=0$ by its stable complement. Then
\[
p,f_p,f_q,q,x,s,y
\]
is induced: labels and boundary stability remove root--cycle and root--root chords, singleton contacts remove root--clique chords, and component separation removes clique--cycle chords.

Choose an unused position $f\in E$ and a partial root $t$ there. It misses $D$, or it belongs to $S$, contradicting the choice of $f$. Each noncompleting $p\in S$ misses $f$. Choose $g\in E-f$ missed by $p$, and a transition $p,x,y$ in $D$ with $py=0$; the matching profile is a nonempty proper subset of connected $D$. If $pt=0$, the path
\[
y,x,p,r,t,f,g
\]
is induced: $r$ misses both zero components, $t$ misses $D$ and $g$, $p$ misses $f,g,y$, and the components are anticomplete. Thus $t$ sees all but at most one vertex of $S$ (the possible completer of $E$), giving at least two neighbors in matching $A$. This excludes $E$, including a zero triangle.

\emph{Four-cliques.} If a zero four-clique $E$ and another zero clique $F$ of order at least three occur, delete the possible completer of $F$. Partial roots still cover $F$ and at least three positions $E'\subseteq E$, since the deleted root sees at most one vertex of $E$. The remaining roots have at most one contact in each of $E',F$, contradicting the crossing rule. Thus $E$ is the only large zero clique; the mixed-component exclusion makes all other zero components bipartite. For distinct $f,g\in E$, prescribe colors $2,1$ on the disjoint stable groups $N_A(f),N_A(g)$. These prescriptions extend to a proper $1,2$-coloring of matching $A$. Give $f,g$ colors $1,2$, the other two clique positions $3,4$, and all other zero components $3,4$. This four-colors $H$, a contradiction.

It remains to color $H$ when every zero component is triangle-free or a triangle; the mixed-component argument also shows that triangles exclude nonbipartite triangle-free components.

\emph{No zero triangles.} Extend $W_0$ to a maximal stable set $W$ of $A$. Color $W$ with $1$ and its stable complement with $2$. Each nonbipartite zero component has at most one color-$2$ boundary root, so the palette rule four-colors $H$.

\emph{Zero triangles occur.} Every other zero component is bipartite, and each triangle position has at least two $A$-neighbors. Let $P$ be all zero-triangle completers. If $d,e\in P$ were adjacent, their disjoint zero profiles would complete distinct triangles $E,F$ and miss the opposite triangle. Deleting $d,e$ leaves partial roots covering both triangles, contradicting the crossing rule. Thus $P$ is stable.

Suppose the mate $q$ of $d\in P$, completing $E$, contacts a zero triangle $F$. Since $P$ is stable and adjacent roots have disjoint zero profiles, $q\notin P$, $q$ misses $E$, and its sole $F$-contact $f$ is missed by $d$. Delete $d$ and the possible completer $c$ of $F$, and put $T=N_A(E)\setminus\{d,c\}$, omitting $c$ if absent. Every $E$ position has a partial root after deletion of $d$; deleting $c$ loses at most one position. Thus $T$ covers at least two positions, and its roots have at most one contact in each of $E,F$. Every $p\in T$ sees $f$, or its contact is crossed with $q$'s, forcing $pq=1$ although $q$'s mate is $d$. Choose $g\in F-f$ missed by $d$. Its at least two $A$-neighbors include a partial root $h\ne d,c$. This root misses $E$: otherwise $h\in T$ sees $f$, contradicting its sole contact $g$. Crossing now makes $h$ adjacent to every member of $T$, impossible in matching $A$ since $|T|\ge2$. Hence the mates $M_P$ of $P$ miss all zero triangles.

Set $A_0=A-(P\cup M_P)$. Its roots are partial or absent on every zero triangle. At most one triangle has all three $A_0$-attachment groups nonempty, by the crossing rule. Choose one group of that triangle, or the empty set if none exists, and extend it to a maximal stable set $W$ of $A_0$. Color $W\cup P$ with $1$ and the rest of $A$ with $2$; this is proper since $P$ is stable, its mates are in $M_P$, and $A_0-W$ is stable. Every triangle has a chosen or empty $A_0$-group; its entire original attachment group has color $1$, since $M_P$ misses triangles. The palette rule colors all triangles with $2,3,4$ and all other zero components with $3,4$. Both cases four-color $H$, a contradiction.
\end{proof}

\subsection{Disjoint maximum four-cliques}

\begin{lemma}
\label{lem:q4-two-triangles-tail}

Suppose $a,r,x_1,x_2,x_3,x_4$ is an induced path in a $\{P_7,\mathrm{diamond}\}$-free graph. The neighborhood of $a$ cannot contain two anticomplete triangles $T,S$ both missed by $r$.

\end{lemma}
\begin{proof}

Every tail vertex has at most one neighbor in each triangle, by the diamond with $a$. Every triangle vertex contacts the tail, or it prefixes the displayed path to give $P_7$. Let $i$ be the first tail position owned by either triangle; then $i\le2$.

Suppose $T$ owns $x_i$ at $b_i$ and $x_k$ at a different $b_k$, where $k\ge i+2$, while $S$ owns $x_k$ at $c_k$ missing $x_i$. Choose $c_0\in S\setminus\{c_k\}$ missing $x_i$. Then

\[
r,x_1,\ldots,x_i,b_i,b_k,x_k,c_k,c_0
\]
is induced and contains $P_7$: earlier tail positions have no owners, the owner conditions remove triangle--tail chords, and the original path and anticomplete triangles remove the rest. If $i=2$, both triangles have distinct owners at positions $2,3,4$, contradicting this path with $k=4$. Thus $i=1$.

For each triangle let $I\subseteq\{3,4\}$ consist of positions whose owner misses $x_1$. Both subsets are nonempty, since all three owners occur, and disjoint by the switching path. A triangle missing $x_1$ would have $I=\{3,4\}$, so $x_1$ has both owners. Exchange the triangles so these subsets are $\{3\}$ and $\{4\}$. Their complete owner patterns are therefore

\[
T:(b_1,b_2,b_3,0\text{ or }b_1),\quad
 S:(c_1,c_2,0\text{ or }c_1,c_3),
\]
with the three named owners distinct in each triangle. All four possibilities give the same induced path $b_3,b_2,x_2,x_1,c_1,c_3,x_4$, a contradiction.
\end{proof}

\begin{lemma}
\label{lem:q4-disjoint-cliques}

Let $G$ be an eight-vertex-critical $\{P_7,\mathrm{diamond}\}$-free graph with $\omega(G)=4$. Then its maximum four-cliques are pairwise disjoint.

\end{lemma}
\begin{proof}

Suppose $a$ is in two maximum four-cliques. The blocks of $N(a)$ have order at most three, and at least two are triangles. Put $D=G-N[a]$. Connectedness gives each component of $D$ a nonempty boundary in $N(a)$.

Exactly one component $E$ of $D$ is not four-colorable. There is at least one, since otherwise four colors on $D$ and $a$, and three new colors on $N(a)$, suffice. Every boundary vertex of such a component is bad, by \cref{lem:q4-good-extension}. If two such components existed, choose boundary vertices $r,s$ and bad tails $r,e_1,e_2,e_3$ and $s,f_1,f_2,f_3$ in them. The two-tail observation in \cref{sec:rooted} makes their boundaries disjoint and makes each miss the other component. If $rs=1$, then $e_2,e_1,r,s,f_1,f_2,f_3$ is induced; if $rs=0$, then $e_3,e_2,e_1,r,a,s,f_1$ is induced. Both contradict $P_7$-freeness.

Fix two triangular blocks $T_1,T_2$ of $N(a)$. Since $d_G(a)\ge7$, choose $r\in N(a)\setminus(T_1\cup T_2)$; the cluster structure makes $r$ miss both. By \cref{lem:q4-two-triangles-tail}, $r$ starts no $P_5$ into $D$. If $r$ contacts $E$, take $w=r$; otherwise take a boundary vertex $w$ of $E$. In the latter case, if $rw=1$, then $w$ misses the same two triangles; if $rw=0$, any rooted $P_5$ from $w$ into $E$ can be prefixed by $r,a$. Thus in every case $w$ contacts $E$ and starts no $P_5$ into it.

Let $L$ be the block of $N(a)$ containing $w$, put $A=N(a)\setminus L$, and consider the triangular blocks of $N_E(w)$. Ownership gives each $A$-vertex at most one contact in such a block. Call a block covered if all its vertices contact $A$. From every block not covered by $A$, choose a vertex missing $A$. If covered blocks occur, fix $b\in A$ contacting one and choose its neighbor in every covered block. These choices exist: \cref{lem:q4-shared-boundaries} confines the $A$-contacts of any two covered blocks to one triangle perfectly matched to both; using the fixed block makes that triangle common to all of them.

The chosen set $X$ is stable and may be empty. At most one vertex of each triangular block $T$ of $A$ contacts $X$: with at most one covered block this follows from ownership of $\{a\}\cup T$; with several, $N_A(X)=\{b\}$. Choose a representative of each such $T$ missing $X$. If $L$ is triangular, choose its representative outside $w$; it misses $N_E(w)$ by the diamond with $a,w$, and hence misses $X$. The representatives and $X$ together are stable.

The set $X$ meets every triangular block of $N_E(w)$, so $N_{E-X}(w)$ is a matching and isolates. Apply \cref{lem:q4-matching-root} with $(C,J,r)=(E,E-X,w)$ to four-color all of $E-X$.

Color $D-X$ with $1,2,3,4$, and color $X$ and all representatives with $5$. The remaining blocks of $N(a)$ have order at most two, so use $6,7$ there and give $a$ color $1$. These palettes respect every boundary and seven-color $G$, a contradiction.
\end{proof}

\subsection{Heavy clique pairs}
\label{sec:q4-three-edge-pair}
\begin{lemma}
\label{lem:q4-three-edge-pair}
\label{lem:q4-perfect-pair}

Let $G$ be connected, $\{P_7,\mathrm{diamond}\}$-free, with $\omega(G)=4$ and $\delta(G)\ge4$. Two disjoint four-cliques with at least three cross edges dominate every vertex belonging to a four-clique: if they are $R,Q$, every four-clique is contained in $N[R\cup Q]$.

\end{lemma}
\begin{proof}

Write the cliques as $R=\{r_i:i\in I\}$ and $Q=\{q_i:i\in I\}$, where $I=\{0,1,2,3\}$ and the cross edges are $r_iq_i$ for $i\in M$, with $M=\{0,1,2\}$ or $I$. Set $K=R\cup Q$, $A=N(K)\setminus K$ and $D=G-N[K]$. A vertex outside $K$ has at most one contact in each clique, by diamond-freeness and $\omega(G)=4$. Its type $(i,j)$ records these contacts in $R,Q$; $\bot$ means no contact and excludes no index.

Every $x\in A$ starts a $P_5$ into $K$: after exchanging sides, its type $(i,j)$ has $i\ne\bot$, and we use $x,r_i,r_a,q_a,q_b$ with $a\in M\setminus\{i,j\}$ and $b\in I\setminus\{i,j,a\}$. A boundary vertex of a component of $D$ is therefore complete to it, or a transition edge prefixes this path to give $P_7$. Connectedness supplies a boundary for each component, so \cref{lem:clique-ownership} makes each component a clique. If a component has an edge, every boundary vertex belongs to the unique maximal clique containing that edge. A component vertex would have degree at most three. Thus $D$ is stable.

Suppose $z\in D$ belongs to a four-clique $C$. Since $\delta(G)\ge4$, choose $y\in N(z)\setminus C$. By \cref{lem:clique-ownership}, $y$ misses $C\setminus\{z\}$, and the three vertices of $C\setminus\{z\}$ have pairwise disjoint core-contact sets. All these neighbors are in $A$, since $D$ is stable. As $y$ has at most two core contacts, some $x\in C\setminus\{z\}$ has contacts disjoint from $y$'s. Write their types as $(i,j),(s,t)$. If the four indices are real and distinct, then $r_j,r_i,x,z,y,q_t,q_s$ is an induced $P_7$.

Otherwise exchange sides or vertices so that $i$ is real and $j\in\{i,\bot\}$ or $t=i$: a one-contact or diagonal type gives the first option; two off-diagonal types give the second. If $a\in M\setminus\{i,j,s,t\}$ exists, choose $b\in I\setminus\{i,j,a,t\}$, possible since at most three real indices are forbidden. Then $y,z,x,r_i,r_a,q_a,q_b$ is induced.

Otherwise the three real indices are exactly $M=\{0,1,2\}$. Reorient so that $i,t$ are real and distinct; otherwise each of $(i,t),(s,j)$ has a missing entry or equal entries, leaving at most two distinct indices. Now $r_3,r_i,x,z,y,q_t,q_3$ is induced. In each case the matching cross-edge sets and disjoint contacts exclude all chords, a contradiction.
\end{proof}

\section{Interfaces between four-cliques}\label{sec:interfaces}

Each local argument fixes a vertex from the three-clique graph and proves that it has ambient degree at most four. \Cref{sec:clique-four} applies these obstructions inside the critical counterexample, then propagates the remaining interfaces through its entire maximum-clique family.

Unless a lemma states otherwise, $A,B,C$ are disjoint four-cliques, and $K=A\cup B\cup C$. An \emph{interface} is the entire edge set between two cliques. We assume every interface is a matching and every vertex outside $K$ has at most one neighbor in each clique. These assumptions hold for disjoint maximum four-cliques in a diamond-free graph of clique number four, by \cref{lem:clique-ownership}. We call such a vertex \emph{exterior}. When the cliques are named $Q,R,T$, the same conventions apply with $K=Q\cup R\cup T$. An interface is \emph{perfect} if it has four edges and \emph{heavy} if it has at least three. A pair dominates the third clique when every vertex of that clique has a neighbor in the pair.

Use $I=\{0,1,2,3\}$, $S=\{1,2,3\}$ and $J=\{2,3\}$ unless local labels are introduced. An absent contact or undefined map value is $\bot$; it is not a real index and excludes no real label. Every optional contact ranges over all labels of its clique and $\bot$. Any target claimed to have no exterior neighbor is fixed from the three-clique graph before considering an exterior vertex or its optional contacts.

The path proofs use three consecutive two-vertex clique blocks and one exterior vertex $z$. For an ordered sequence $P$, let $X(P)$ be the complete edge set between its selected vertices in different cliques, and let $Z(P)$ be its selected neighbors of $z$. The sequence is an induced $P_7$ precisely when these are its consecutive cross-clique edges and the vertices consecutive to $z$, respectively: the three same-clique pairs already give the other three consecutive edges. Thus a complete $X(P),Z(P)$ listing also specifies every required nonedge. \cref{local-crossed-contacts} states its own more general pair-relative hypotheses.

\subsection{Two perfect interfaces}

\begin{lemma}
\label{local-two-perfect}
If $A$--$B$ and $B$--$C$ are perfect, every vertex of $A$ missed by $C$ has no exterior neighbor in a $P_7$-free ambient graph, and has degree four.

\end{lemma}
\begin{proof}
Normalize $E(A,B),E(B,C)$ as $a_ib_i,b_ic_i$ for $i\in I$, and fix the given vertex $a_0$. Suppose an exterior neighbor $z$ exists. Write $E(A,C)=\{a_ic_{\pi(i)}:i\in\operatorname{dom}\pi\}$, where $\pi$ is any partial injection with $\pi(0)=\bot$. Let $j,k\in I\cup\{\bot\}$ be the optional $B,C$ contacts of $z$.

We use three templates:
\[
\begin{aligned}
E(p,q,r)&=(z,a_0,a_p,b_p,b_q,c_q,c_r),\\
M(p,q;h)&=(b_q,b_p,a_p,a_0,z,c_k,c_\ell),\\
R(u,v)&=(z,b_j,b_0,c_0,c_v,a_u,a_v).
\end{aligned}
\]
Their sufficient conditions and complete selected edge sets are as follows.

\begin{enumerate}[label=(\roman*),leftmargin=2.2em,itemsep=0.25em,topsep=0.3em]
\item If $p\in S\setminus\{j\}$, $q\in S\setminus\{j,\allowbreak p,\allowbreak k,\allowbreak \pi(p)\}$, and $r\in I\setminus\{p,\allowbreak q,\allowbreak k,\allowbreak \pi(p)\}$, take $P=E(p,\allowbreak q,\allowbreak r)$; then $X(P)=\{a_pb_p,\allowbreak b_qc_q\}$ and $Z(P)=\{a_0\}$.
\item If $\{p,\allowbreak q,\allowbreak h\}=S$, $j\notin\{p,\allowbreak q\}$, $\{k,\allowbreak \ell\}=\{0,\allowbreak h\}$, and $\pi(p)\notin\{0,\allowbreak h\}$, take $P=M(p,\allowbreak q;h)$; then $X(P)=\{a_pb_p\}$ and $Z(P)=\{a_0,\allowbreak c_k\}$.
\item If $j\in S$, $\{u,\allowbreak v\}=S\setminus\{j\}$, $\pi(u)=v$, $\pi(v)\ne0$, and $k\notin\{0,\allowbreak v\}$, take $P=R(u,\allowbreak v)$; then $X(P)=\{b_0c_0,\allowbreak a_uc_v\}$ and $Z(P)=\{b_j\}$.
\end{enumerate}

The identity matchings give exactly the displayed $A$--$B$ and $B$--$C$ edges. For $E,M$, the exclusions of $\pi(p)$ and $\pi(0)=\bot$ remove every $A$--$C$ edge. For $R$, injectivity gives $\pi(v)\ne v$, so its sole selected $A$--$C$ edge is $a_uc_v$. The contact exclusions give exactly $Z(P)$, including absent contacts. Thus each template satisfies the full selected-edge criterion.

If $k$ is real and $j,k$ are not distinct members of $S$, choose $h\in S$ with $j\in\{\bot,0,h\}$ and $k\in\{0,h\}$. Put $\{p,q\}=S\setminus\{h\}$ and $\{k,\ell\}=\{0,h\}$. If either $\pi(p)$ or $\pi(q)$ lies outside $\{0,h\}$, orient it as $p$ and use $M(p,q;h)$. Otherwise injectivity gives $\{\pi(p),\pi(q)\}=\{0,h\}$; orient so $\pi(p)=k$ and use $E(p,q,\ell)$. This includes every $k=0$ profile.

Suppose $k=\bot$. For $j\in\{0,\bot\}$ choose $p\in S$ with $\pi(p)\ne0$, then $q\in S\setminus\{p,\pi(p)\}$, and use $E(p,q,0)$. For $j\in S$, put $\{p,q\}=S\setminus\{j\}$. Use $E(p,q,0)$ if $\pi(p)\notin\{0,q\}$, interchanging $p,q$ if possible. Otherwise $\pi(p)\in\{0,q\}$ and $\pi(q)\in\{0,p\}$. If one value is zero, orient so $\pi(p)=0$ and use $D$ below; if neither is zero, use $R(q,p)$.

Finally let $j,k\in S$ be distinct, and let $\{t\}=S\setminus\{j,k\}$. If $\pi(k)\ne t$, choose $r\in\{0,j\}\setminus\{\pi(k)\}$ and use $E(k,t,r)$. Otherwise use $R(k,t)$ when $\pi(t)\ne0$, and $G$ when $\pi(t)=0$:

For $P=D=(a_p,\allowbreak a_0,\allowbreak z,\allowbreak b_j,\allowbreak b_q,\allowbreak c_q,\allowbreak c_p)$, we have $X(P)=\{b_qc_q\}$ and $Z(P)=\{a_0,\allowbreak b_j\}$.

For $P=G=(c_t,\allowbreak c_0,\allowbreak a_t,\allowbreak a_0,\allowbreak z,\allowbreak b_j,\allowbreak b_k)$, we have $X(P)=\{a_tc_0\}$ and $Z(P)=\{a_0,\allowbreak b_j\}$.

For $D$, $\pi(p)=0$ removes all selected $A$--$C$ edges; for $G$, $\pi(t)=0$ gives exactly $a_tc_0$. The distinct indices, $\pi(0)=\bot$, and the specified $j,k$ give every other entry of $X(P),Z(P)$. Injectivity supplies $p$ in the first $k=\bot$ case; the choices of $q,r$ remove at most two of three or one of two labels. All displayed paths have seven distinct vertices. Every profile therefore gives an induced $P_7$. The fixed $a_0$ has no exterior neighbor, leaving just its three $A$-neighbors and $b_0$, so its ambient degree is four.
\end{proof}

\subsection{Two heavy interfaces with a common clique}

\begin{lemma}
\label{local-common-missing-label}
\label{local-common-central-defect}
\label{local-one-three-three}
\label{local-aligned-defects}
\label{local-two-three-three}
\label{local-sparse-outer}
Under the local interface conventions, suppose $A$--$B$ and $A$--$C$ each have three edges, $B$--$C$ has at most two edges, and each heavy pair dominates the remaining clique. In a $P_7$-free ambient graph, the vertex of $A$ missed by $C$ has no exterior neighbor. Its ambient degree is three if the two unmatched $A$-vertices coincide, and four otherwise.

\end{lemma}
\begin{proof}
Domination makes the outer interface nonempty. The target is fixed from the two heavy matchings before any exterior profile. First suppose the two unmatched vertices of $A$ coincide at $a_0$. Normalize $E(A,B)=\{a_ib_i:i\in S\}$ and $E(A,C)=\{a_ic_i:i\in S\}$. Domination makes $b_0,c_0$ incident with the nonempty $B$--$C$ interface. Fix an exterior neighbor $z$ of $a_0$; ownership makes $z$ miss every other $A$-vertex.

If $z$ sees $b_s$, let $t\in I\cup\{\bot\}$ be its optional $C$-contact. For every $i\in S\setminus\{s,t\}$, \cref{local-crossed-contacts}, with center $a_0$ and roots $z,a_i$, forces $b_sc_i=1$, since $za_i=0$. A matching allows at most one such $i$, so $s,t$ must be distinct members of $S$. For the remaining $i$, the same lemma with roots $a_i,z$ forces $b_ic_t=1$. These two edges exhaust the interface but miss $b_0,c_0$, contradicting domination. Interchanging $B,C$ covers a $C$-contact without a $B$-contact. Thus $z$ misses $B\cup C$.

If the interface has one edge, it is $b_0c_0$. For distinct $i,k\in S$, the sequence $z,a_0,a_i,b_i,b_0,c_0,c_k$ is induced: its complete core cross-edge set is $\{a_ib_i,b_0c_0\}$ and $z$ sees only $a_0$. Otherwise there are two interface edges. Choose an interface edge $b_jc_k$ with $k\in S$, possible because there are two matching edges. Choose $i\in S\setminus\{j,k\}$ and
\[
\ell\in I\setminus\bigl(\{i,j\}\cup\{h:b_hc_i=1\}\bigr).
\]
At most two of three indices, and three of four, respectively, are excluded. Now $z,\allowbreak a_0,\allowbreak a_i,\allowbreak c_i,\allowbreak c_k,\allowbreak b_j,\allowbreak b_\ell$ is induced. Its only core cross edges are $a_ic_i,c_kb_j$: the chosen indices remove the other $A$-contacts, and matching plus the choice of $\ell$ removes every other selected $B$--$C$ edge. The root sees only $a_0$, so the selected-edge criterion gives the contradiction.

Thus the core-selected $a_0$ has no exterior neighbor and has ambient degree three.

Now suppose the unmatched vertices of $A$ differ. Within this subcase put
\[
A=\{a_1,a_2,a_3,a_4\},\quad B=\{b_0,b_1,b_2,b_3\},\quad
 C=\{c_0,c_1,c_2,c_4\},
\]
with $A$--$B$ edges $a_ib_i$ for $i=1,2,3$ and $A$--$C$ edges $a_ic_i$ for $i=1,2,4$. Domination makes $b_0,c_0$ incident with the nonempty interface of size at most two. Fix an exterior neighbor $z$ of $a_3$. Let $B^- =\{b_0,b_1,b_2\}$ and let $b_s,c_t$ be its optional contacts in $B^-,C$. Here $s\in\{\bot,0,1,2\}$ and $t\in\{\bot,0,1,2,4\}$. The value $s=\bot$ allows a contact with $b_3$; all sequences in that case omit $b_3$.

\emph{No contact in $B^-$.} If $t$ is real, apply \cref{local-crossed-contacts} to center $a_3$, cliques $B^-,C$, and roots $a_i,z$ for $i\in\{1,2\}\setminus\{t\}$. It forces $b_ic_t=1$. Matching therefore gives $t\in\{1,2\}$, and, with $i$ its other member, the interface is $\{b_ic_t,b_0c_0\}$. Use
\[
a_4,a_3,z,c_t,c_0,b_0,b_t.
\]
Its only core cross edge is $c_0b_0$, and $z$ sees exactly $a_3,c_t$. If $t=\bot$, let $b_0c_d$ be the forced interface edge. For $d\ne4$, take $i=4$ and choose $\ell\in\{1,2\}$ with $b_\ell c_4=0$. For $d=4$, take $i=1,\ell=2$; then $c_1$ misses $B$, since the other interface edge must meet $c_0$. In either case use
\[
z,a_3,a_i,c_i,c_d,b_0,b_\ell.
\]
Its only core cross edges are $a_ic_i,c_db_0$, and $z$ sees only $a_3$. The choices give distinct within-clique indices and exclude the other cross edges. In particular $b_3$ is absent from both sequences.

\emph{A real contact $b_s$ and $t\ne4$.} Use \cref{local-crossed-contacts} on $B^-,C$ with center $a_3$ and roots $z,a_4$; the latter misses $B^-$. It forces $b_sc_4=1$. For $s=h\in\{1,2\}$, let $k$ be its other member. The interface is $\{b_hc_4,b_0c_0\}$. If $t=0$, use the first sequence below; otherwise use the second:
\[
c_h,c_k,a_k,a_3,z,b_h,b_0;\quad
 z,a_3,a_4,c_4,c_0,b_0,b_k.
\]
Their core cross edges are respectively $c_ka_k$ and $a_4c_4,c_0b_0$; the selected neighbors of $z$ are $a_3,b_h$ and $a_3$. All other selected cross pairs are nonedges.

For $s=0$, the interface is $\{b_0c_4,b_pc_0\}$, $p\in\{1,2,3\}$. If $p=3$, choose $i\in\{1,2\}\setminus\{t\}$, let $j$ be its other member, choose $k\in\{0,j\}\setminus\{t\}$, and use
\[
b_j,b_0,z,a_3,a_i,c_i,c_k.
\]
Its only core cross edge is $a_ic_i$; $z$ sees $b_0,a_3$. Each two-element choice excludes at most one index. For $p\in\{1,2\}$ and $t=0$, use
\[
a_4,a_1,c_1,c_0,z,b_0,b_3.
\]
Its only core cross edge is $a_1c_1$, and $z$ sees $c_0,b_0$. For $p\in\{1,2\}$ and $t\ne0$, choose $\ell\in\{1,2\}\setminus\{t\}$ and use
\[
a_4,a_3,z,b_0,b_p,c_0,c_\ell.
\]
The only core cross edge is $b_pc_0$, and $z$ sees $a_3,b_0$. The displayed interface and the omitted indices exclude all other cross edges in these sequences, including the possible $a_3b_3$ edge.

\emph{A real contact $b_s$ and $t=4$.} Choose $i\in\{1,2\}\setminus\{s\}$. If $b_ic_4=1$, domination gives the interface $\{b_ic_4,b_0c_0\}$. Otherwise put
\[
X=B^-\setminus(\{b_i,b_s\}\cup N_B(c_4)),\quad
 Y=C\setminus(\{c_4,c_i\}\cup N_C(b_i)).
\]
Here $Y\ne\varnothing$. A nonedge $xy$ with $x\in X,y\in Y$ would give the induced sequence $x,b_i,a_i,a_3,z,c_4,y$: its only core cross edge is $b_ia_i$, and the profiles and set exclusions remove every other cross edge and unwanted $z$-contact. If $X=\varnothing$, $c_4$ has its $B$-neighbor in $B^-\setminus\{b_i,b_s\}$. If $X\ne\varnothing$, matching and the absence of such a nonedge force $X,Y$ to be singletons. Then $b_i$'s interface edge and the $X$--$Y$ edge exhaust the interface; neither meets $c_4$. Consequently $X=B^-\setminus\{b_i,b_s\}$, and domination forces $X=\{b_0\}$ and the two $C$-endpoints to be $c_0,c_s$. Thus all possibilities are the following entire interfaces. In the last two rows $\{h,k\}=\{1,2\}$ and $p\in\{1,2,3\}$:
\[
\begin{array}{c|l|l}
 s& E(B,C)&\text{induced sequence}\\\hline
 0&\{b_hc_4,b_0c_0\},\ h\in\{1,2\}
   &a_1,a_2,c_2,c_4,z,b_0,b_3\\
 h&\{b_kc_4,b_0c_0\}\text{ or }\{b_kc_h,b_0c_0\}
   &c_h,c_k,a_k,a_3,z,b_h,b_0\\
 h&\{b_0c_4,b_pc_0\}\text{ or }\{b_kc_0,b_0c_h\}
   &z,c_4,c_h,a_h,a_k,b_k,b_3.
\end{array}
\]
The core cross edges in these rows are respectively $a_2c_2$, $c_ka_k$, and $c_ha_h,a_kb_k$. All selected $B$--$C$ pairs are nonedges by the displayed entire interfaces. The selected neighbors of $z$ are respectively $c_4,b_0$; $a_3,b_h$; and $c_4$. The remaining matching nonedges and distinct indices follow from $\{h,k\}=\{1,2\}$. Every sequence above therefore satisfies the selected-edge criterion. The cases exhaust both real and absent contacts, including an allowed $b_3$ contact when $s=\bot$, and each gives the required contradiction.

The deductions used $|E(B,C)|\le2$, not equality. For the sole edge $b_0c_0$, the no-contact branch has $d=0$; a real $B^-$-contact with $t\ne4$ forces a second edge, while $t=4$ leaves candidate sets of sizes at least one and two, hence a nonedge sequence. Thus the one-edge case is included. The fixed $a_3$ has no exterior neighbor, and its three $A$-mates and $b_3$ give ambient degree four.
\end{proof}

\subsection{A single heavy interface}

\begin{lemma}
\label{local-single-heavy}
Suppose $Q$--$R$ is heavy, $Q$--$T$ and $R$--$T$ have at most two edges each, and $Q\cup R$ dominates $T$. In a $P_7$-free ambient graph, a vertex of $T$ chosen from the three-clique graph has degree four.

\end{lemma}
\begin{proof}
Both smaller interfaces have size two, and each $T$-vertex has exactly one neighbor in $Q\cup R$. Label the heavy matching as $q_ir_i$ for $i\in I\setminus\{d\}$, where $d=\bot$ denotes the perfect case. This is a labeling of the actual edges; $q_dr_d$ remains absent when $d$ is real. Let $L,H\subset I$ be the two-element sets of owner indices in $Q,R$, respectively. Write $a_i\in T$ for the vertex owned by $q_i$, and $b_j\in T$ for the vertex owned by $r_j$.

First suppose $L\cap H\ne\varnothing$. Label an intersection index $0$, choose the target $a_0$, and let $b_k$ be the other $R$-owned $T$-vertex, so $k\in S$. For an exterior neighbor $z$ of $a_0$ let $s,t$ be its optional $Q,R$-contact indices. All possibilities are as follows.

\begin{enumerate}[label=(\roman*),leftmargin=2.2em,itemsep=0.25em,topsep=0.3em]
\item $t\ne0$. Choose $i\in S\setminus\{s,\allowbreak t,\allowbreak d\}$ and $j\in S\setminus\{s,\allowbreak i\}$. Take $P=(z,\allowbreak a_0,\allowbreak b_0,\allowbreak r_0,\allowbreak r_i,\allowbreak q_i,\allowbreak q_j)$; then $X(P)=\{b_0r_0,\allowbreak r_iq_i\}$ and $Z(P)=\{a_0\}$.
\item The first choice of $i$ fails: $s,\allowbreak t,\allowbreak d$ are distinct and exhaust $S$. Take $P=(q_d,\allowbreak q_s,\allowbreak z,\allowbreak a_0,\allowbreak b_0,\allowbreak r_0,\allowbreak r_d)$; then $X(P)=\{b_0r_0\}$ and $Z(P)=\{q_s,\allowbreak a_0\}$.
\item $t=0$. Choose $i\in S\setminus\{s,\allowbreak k,\allowbreak d\}$ and $j\in S\setminus\{s,\allowbreak i\}$. Take $P=(q_j,\allowbreak q_i,\allowbreak r_i,\allowbreak r_0,\allowbreak z,\allowbreak a_0,\allowbreak b_k)$; then $X(P)=\{q_ir_i\}$ and $Z(P)=\{r_0,\allowbreak a_0\}$.
\item This choice of $i$ fails: $s,\allowbreak k,\allowbreak d$ are distinct and exhaust $S$. Take $P=(q_d,\allowbreak q_s,\allowbreak z,\allowbreak a_0,\allowbreak b_k,\allowbreak r_k,\allowbreak r_d)$; then $X(P)=\{b_kr_k\}$ and $Z(P)=\{q_s,\allowbreak a_0\}$.
\end{enumerate}

Whenever $i$ exists, so does $j$, since at most two members of $S$ are excluded. In the exceptional cases $d$ is real and $q_dr_d$ is absent; every other selected $Q,R$ pair has unequal indices. In all four cases $q_0$, the owner of $a_0$, is absent, and the remaining $T$-owners and contact exclusions give exactly the two specified edge sets.

Now suppose $L\cap H=\varnothing$. Exchange $Q,R$ if necessary, and label so $L=\{0,1\}$, $H=\{2,3\}$, and the only possibly missing heavy edge is $q_1r_1$. Fix the target $a_0$. For an exterior neighbor with contact indices $s,t$ use:

\begin{enumerate}[label=(\roman*),leftmargin=2.2em,itemsep=0.25em,topsep=0.3em]
\item $s\ne1$ and $\{s,\allowbreak t\}\ne J$. Choose $i\in J\setminus\{s,\allowbreak t\}$ and $j\in\{0,\allowbreak 2,\allowbreak 3\}\setminus\{i,\allowbreak t\}$. Take $P=(z,\allowbreak a_0,\allowbreak a_1,\allowbreak q_1,\allowbreak q_i,\allowbreak r_i,\allowbreak r_j)$; then $X(P)=\{a_1q_1,\allowbreak q_ir_i\}$ and $Z(P)=\{a_0\}$.
\item $s=1,\allowbreak t\ne0$. Choose $k\in J\setminus\{t\}$, $\{i\}=J\setminus\{k\}$. Take $P=(q_i,\allowbreak q_1,\allowbreak z,\allowbreak a_0,\allowbreak b_k,\allowbreak r_k,\allowbreak r_0)$; then $X(P)=\{b_kr_k\}$ and $Z(P)=\{q_1,\allowbreak a_0\}$.
\item $s=1,\allowbreak t=0$. Write $\{j,\allowbreak k\}=J$. Take $P=(q_k,\allowbreak q_1,\allowbreak z,\allowbreak r_0,\allowbreak r_j,\allowbreak b_j,\allowbreak b_k)$; then $X(P)=\{r_jb_j\}$ and $Z(P)=\{q_1,\allowbreak r_0\}$.
\item $\{s,\allowbreak t\}=J$. Take $P=(z,\allowbreak r_t,\allowbreak r_0,\allowbreak q_0,\allowbreak q_1,\allowbreak a_1,\allowbreak b_s)$; then $X(P)=\{r_0q_0,\allowbreak q_1a_1\}$ and $Z(P)=\{r_t\}$.
\end{enumerate}

In the first case $i$ fails to exist only for the excluded pair $\{s,t\}=J$, and $j$ always exists. The other choices remove at most one of two indices. No case selects $r_1$; thus the possible missing edge is never used or added. The matching and single-owner formulas give exactly $X(P)$, while $s,t$ give exactly $Z(P)$, including the last two cases where the target $a_0$ is absent from the path. Every profile has produced an induced $P_7$. Consequently the fixed target has no exterior neighbor, and its three $T$-neighbors and one owner give degree four.
\end{proof}

\subsection{A perfect and a three-edge interface}

\begin{lemma}
\label{local-perfect-three}
Suppose
\[
E(A,B)=\{a_ib_i:i\in I\},\quad
E(B,C)=\{b_ic_i:i\in S\},
\]
where $I=\{0,1,2,3\}$ and $S=\{1,2,3\}$, and $A\cup B$ dominates $C$. In a $P_7$-free ambient graph, $c_0$ has no exterior neighbor and has degree four, for every matching $A$--$C$.

\end{lemma}
\begin{proof}
Let $\sigma(p)$ be the $A$-owner index of $c_p$, or $\bot$ if absent. It is a partial injection, with $\sigma(0)=x\in I$ by domination. If $A$--$C$ is perfect, apply \cref{local-two-perfect} to $(C,A,B)$. Otherwise fix any $d\in S$ with $\sigma(d)=\bot$. Both $c_0$ and $d$ are fixed before an exterior neighbor $z$ of $c_0$ is considered. Let $s,t\in I\cup\{\bot\}$ be its optional $A,B$ contacts.

The following templates satisfy the full selected-edge criterion:

\begin{enumerate}[label=(\roman*),leftmargin=2.2em,itemsep=0.25em,topsep=0.3em]
\item Choose $w\in S\setminus\{t\}$ and distinct $q,\allowbreak r\in I\setminus\{w,\allowbreak x,\allowbreak s,\allowbreak \sigma(w)\}$ with $q\ne t$. Take $P=E(w;q,\allowbreak r)=\allowbreak (z,\allowbreak c_0,\allowbreak c_w,\allowbreak b_w,\allowbreak b_q,\allowbreak a_q,\allowbreak a_r)$; then $X(P)=\{c_wb_w,\allowbreak b_qa_q\}$ and $Z(P)=\{c_0\}$.
\item For $s\notin\{\bot,\allowbreak x\}$, choose distinct $q,\allowbreak r\in I\setminus\{s,\allowbreak d,\allowbreak t\}$ with $q\ne x$. Take $P=Y(q,\allowbreak r)=\allowbreak (b_r,\allowbreak b_q,\allowbreak a_q,\allowbreak a_s,\allowbreak z,\allowbreak c_0,\allowbreak c_d)$; then $X(P)=\{b_qa_q\}$ and $Z(P)=\{a_s,\allowbreak c_0\}$.
\item For $t\in I$, choose $w\in S\setminus\{t\}$ and distinct $u,\allowbreak v\in I\setminus\{t,\allowbreak x,\allowbreak \sigma(w),\allowbreak s\}$ with $v\ne w$. Take $P=V(u,\allowbreak v,\allowbreak w)=\allowbreak (a_u,\allowbreak a_v,\allowbreak b_v,\allowbreak b_t,\allowbreak z,\allowbreak c_0,\allowbreak c_w)$; then $X(P)=\{a_vb_v\}$ and $Z(P)=\{b_t,\allowbreak c_0\}$.
\end{enumerate}

For $E,V$, the selected $A$-indices avoid $x,\sigma(w)$; for $Y$ they avoid $x$, and $c_d$ is unowned. These remove all $A$--$C$ edges. The other index exclusions give exactly the identity-matching edges and $s,t$ contacts displayed above, including absent contacts.

First let $s\notin\{\bot,x\}$. If $s,d,t$ are not three distinct real indices, $I\setminus\{s,d,t\}$ has at least two members; choose $q\ne x$ and another $r$ there for $Y$. Otherwise let $j$ be the fourth index. For $x=t$ use $V(d,j,d)$; for $x\ne t$ use $F$ below. Thus every real restricted $A$-contact is excluded, without inspecting any further value of $\sigma$.

Now let $s\in\{\bot,x\}$. If $t\ne d$, use $E(d;q,r)$ with two indices outside $\{d,x\}$, oriented so $q\ne t$. For $t=d$, if some $w\in S\setminus\{d\}$ leaves at least two indices in $I\setminus\{d,x,\sigma(w)\}$, choose $u,v$ there, orienting $v\ne w$, and use $V$. Otherwise $x\ne d$ and $\sigma$ maps $S\setminus\{d\}$ bijectively onto $L=I\setminus\{x,d\}$, by injectivity. If $x\ne0$, take $w=x$; if $x=0$ and some $w\in L$ has $\sigma(w)=w$, take it. In either case $I\setminus\{w,x,\sigma(w)\}$ has two members, so use $E(w;q,r)$, oriented so $q\ne d$. Only $x=0$ with $\sigma$ interchanging the two members $u,v$ of $L$ remains; use $H$.

\begin{enumerate}[label=(\roman*),leftmargin=2.2em,itemsep=0.25em,topsep=0.3em]
\item If $s,\allowbreak d,\allowbreak t$ are distinct real indices, $j\in I\setminus\{s,\allowbreak d,\allowbreak t\}$, and $x\notin\{s,\allowbreak t\}$, take $P=F=\allowbreak (a_t,\allowbreak a_s,\allowbreak z,\allowbreak c_0,\allowbreak c_d,\allowbreak b_d,\allowbreak b_j)$; then $X(P)=\{c_db_d\}$ and $Z(P)=\{a_s,\allowbreak c_0\}$.
\item If $x=0$, $s\in\{\bot,\allowbreak 0\}$, $t=d$, $L=\{u,\allowbreak v\}$, $\sigma(u)=v$, and $\sigma(v)=u$, take $P=H=\allowbreak (a_u,\allowbreak a_v,\allowbreak c_u,\allowbreak c_0,\allowbreak z,\allowbreak b_d,\allowbreak b_0)$; then $X(P)=\{a_vc_u\}$ and $Z(P)=\{c_0,\allowbreak b_d\}$.
\end{enumerate}

For $F$, neither selected $A$-index equals $x$, and $c_d$ is unowned; the $A,B$ index sets are disjoint. For $H$, $u,v,d$ are the three nonzero indices: the only $A$--$C$ edge is $a_vc_u$, and the $A,B$ index sets are disjoint. Its $B$--$C$ pairs are nonedges because $d\ne u$ and $b_0,c_0$ are unmatched across that interface. The specified contacts give exactly $Z(P)$. All sequences have seven distinct vertices, so every profile yields an induced $P_7$. The fixed $c_0$ therefore has no exterior neighbor; its three $C$-neighbors and owner $a_x$ give degree four.
\end{proof}

\section{Maximum four-cliques}\label{sec:clique-four}

\begin{proposition}
\label{lem:small-family}
\label{lem:three-family}
\label{lem:four-perfect-no-triangle}
\label{cor:perfect-clique-family-transversal}
\label{lem:uniform-patterns}
\label{prop:q4-exclusion}

Let $G$ be connected and $\{P_7,\mathrm{diamond}\}$-free, with $\omega(G)=4$, $\delta(G)\ge5$, and pairwise disjoint maximum cliques. If some two maximum cliques have at least three cross edges, there is a stable set meeting every maximum clique.

\end{proposition}
\begin{proof}
Let $\mathcal Q$ be the family of all maximum four-cliques. Ownership makes their interfaces matchings. If $|\mathcal Q|\le4$, choose a stable transversal greedily: each earlier representative forbids at most one position in the next clique. Assume henceforth that $|\mathcal Q|\ge5$. By hypothesis a heavy pair exists, and every heavy pair dominates every third member by \cref{lem:q4-three-edge-pair}.

\emph{Uniform triples.} Every triple with a heavy interface has sizes $(3,3,3)$ or $(4,4,4)$. Ownership makes all interfaces matchings. Two perfect interfaces and one nonperfect interface contradict \cref{local-two-perfect}, which supplies a missing endpoint of ambient degree four. An empty interface is also impossible: the heavy pair is incident with it, so domination forces the other interface to be perfect; domination by this new heavy pair then forces the remaining incident interface to be perfect, giving the preceding contradiction.

Thus every interface is nonempty. If there is exactly one heavy interface, the other two have size at most two and their union dominates the four vertices of the remaining clique. Both consequently have size two, and \cref{local-single-heavy} supplies a vertex of ambient degree four, a contradiction.

Suppose there are at least two heavy interfaces. If a perfect interface occurs and the triple is not all perfect, a three-edge interface also occurs, by the first paragraph. Then \cref{local-perfect-three} applies, again giving ambient degree four. If none is perfect, choose two three-edge interfaces. When the remaining interface has at most two edges, \cref{local-common-missing-label} gives a vertex of ambient degree at most four. The only remaining possibility is $(3,3,3)$, completing the alternative.

\emph{Propagation.} Fix a heavy pair $Q_0,Q_1$ with interface size $s\in\{3,4\}$. Apply the triple conclusion to $Q_0,Q_1,R$ to make every interface incident with either base clique have size $s$. Apply it again to $Q_0,R,S$ to make every remaining interface have size $s$.

\emph{Three-edge interfaces.} If $s=3$, each $R\ne Q_0$ leaves one vertex $u_R$ of $Q_0$ unowned. Distinct $R,S$ have $u_R\ne u_S$, since $R\cup S$ dominates $Q_0$. There are at most four such members, so the family has exactly five. Choose the first representative of $Q_0$ unowned in the fifth member, then choose representatives in the second, third and fourth greedily. At most three vertices of the fifth member are forbidden. This is the required stable transversal.

\emph{Perfect interfaces.} It remains that all interfaces are perfect. We show that no triangle meets three family members. Suppose a triangle meets $Q_0,Q_1,Q_2\in\mathcal Q$, and choose a fourth clique $Q_3\in\mathcal Q$. Label each clique by $0,1,2,3$, normalizing all interfaces from $Q_0$ to the identity. Write $p,q,r$ for the interface permutations from $Q_1$ to $Q_2$, $Q_1$ to $Q_3$, and $Q_2$ to $Q_3$, respectively. A vertex is denoted $(i,a)$ when it has label $a$ in $Q_i$. The triangle gives a fixed point $c$ of $p$.

For every such fixed point, $q(c)\ne c$, $r(c)\ne c$, and $q(c)\ne r(c)$. Indeed, an edge of the triangle cannot have a second common neighbor: two nonadjacent common neighbors give a diamond, and adjacent ones give a four-clique meeting two of the $Q_i$. Put $d=r^{-1}(c)$ and $e=q^{-1}(c)$. Thus $r(d)=q(e)=c$ and $d,e\ne c$. If $b$ is one of the two labels outside $\{c,d\}$ and $p(b)\ne d$, $q(b)\ne c$, let $a$ be the other label. Then

\[
(0,a),(0,b),(1,b),(1,c),(2,c),(2,d),(3,c)
\]
is an induced $P_7$: its possible chords are excluded by the distinct labels, the matching interfaces, and the displayed inequalities. Consequently

\begin{equation}\label{eq:four-perfect-condition}
\{p^{-1}(d),e\}=\{0,1,2,3\}\setminus\{c,d\}.
\end{equation}
In particular $c,d,e$ are distinct. Let $f$ be the remaining label. Then $p(f)=d$. Interchanging $Q_1,Q_2$ replaces $(p,d,e)$ by $(p^{-1},e,d)$, so the same equation gives $p(e)=f$.

If $q(c)=e$, then

\[
(1,e),(1,c),(3,e),(3,f),(0,f),(0,d),(2,d)
\]
is induced. If $q(c)=d$, put $x=r^{-1}(d)$ and let $y$ be the other member of $\{e,f\}$; here $x\in\{e,f\}$ because $r(d)=c$ and $r(c)\ne q(c)$. Then

\[
(0,x),(0,y),(2,y),(2,c),(1,c),(3,d),(3,c)
\]
is induced. In each path, the matching interfaces and established identities exclude every chord. Hence $q(c)=f$. Interchanging $Q_1,Q_2$ leaves $f$ unchanged and gives $r(c)=f$, contrary to $q(c)\ne r(c)$.

Choose a base clique and normalize all its interfaces to the identity. Every other interface is fixed-point-free, since a fixed point would give the excluded triangle. Choose distinct labels $a,b$: the base vertex of label $a$, together with label $b$ in every other member, is a stable transversal. Its base-to-other pairs have unequal labels, and its other pairs have the same label and hence are nonedges. This proves the proposition.
\end{proof}

\begin{corollary}
\label{cor:seven-critical-reduction}
If $\chi\le\omega+3$ fails for $\{P_7,\mathrm{diamond}\}$-free graphs, it fails on a seven-vertex-critical graph of clique number three.

\end{corollary}
\begin{proof}
\cref{lem:minimal-counterexample} supplies a connected $(q+4)$-vertex-critical counterexample with $q\le4$, minimum degree at least $q+3$, and no stable maximum-clique transversal. If $q=4$, \cref{lem:q4-disjoint-cliques} gives disjoint maximum cliques, and \cref{lem:minimal-counterexample} gives a heavy pair. \cref{lem:small-family} then supplies the forbidden transversal. Clique number at most two is excluded by \cref{lem:triangle-free-five} and the edgeless convention. Thus $q=3$ and the critical chromatic number is seven.
\end{proof}

\section{Clique number three}\label{sec:clique-three}

Except in the critical-profile toolkit below, let $G$ be $7$-vertex-critical, $\{P_7,\mathrm{diamond}\}$-free, with $\omega(G)=3$. Thus $G$ is connected, has no clique cutset, and $\delta(G)\ge6$. By \emph{triangle sparsity} we mean \cref{lem:clique-ownership}'s consequences: neighborhoods are matchings and isolates, exterior vertices have at most one contact on a triangle, and edges have at most one common neighbor.

For any triangle $Q=\{a,b,c\}$ put
\[
A=N(a)\setminus Q,\quad B=N(b)\setminus Q,\quad C=N(c)\setminus Q,\quad Z=V(G)\setminus N[Q].
\]
These are the full owner classes and zero set. The three owner-plus-zero domains are $G[A\cup Z]$, $G[B\cup Z]$, and $G[C\cup Z]$. A root in $B$ starts no $P_5$ into $A\cup Z$: prefix such a path by $c,b$, whose only later contact is $b$ with the root. Permuting the owners gives all corresponding no-$P_5$ statements, including those into $Z$.

We first three-color these domains. We then use an owner odd cycle to replace $Q$ by a triangle whose global zero set is smaller.

\hypertarget{rule:two-boundary-colors}{}
\emph{Two-boundary-color palette.} Suppose an induced graph $H$ has full boundary $T_1\cup T_2$, with nonempty classes already colored $1,2$, respectively. If every boundary profile is a vertex cover of $H$, and the vertices of $H$ contacting both classes form a stable set, extend the coloring by giving each vertex a class color whose entire class it misses, or $3$ if it contacts both. Each of the first two color classes is stable because any root in the corresponding boundary class covers every edge; the third is stable by hypothesis. The assignment also respects every boundary edge.

\subsection{Full boundaries}

\begin{lemma}
\label{lem:bipartite-boundary}
\label{lem:owner-domain}
\label{lem:zero-core}
Let $D$ be connected and induced, with at least two vertices, and $S=N_G(D)\setminus D$. Every component of $G-N_G[D]$ is four-colorable. Moreover $G$ is six-colorable if either

\begin{enumerate}[label=(\arabic*)]
\item $\chi(D)\le4$ and $S$ is bipartite; or
\item $\chi(D)\le3$ and the good and bad parts of $S$ on $D$ are, in either order, stable and bipartite.

\end{enumerate}
Every owner domain $G[A\cup Z]$ is four-colorable. If $D$ is a four-chromatic component of it and $J\subseteq D$ is induced and $4$-vertex-critical, then $B_D=N_B(D)$ and $C_D=N_C(D)$ have nonbipartite union. Every root in this union contacts this same $J$ and starts no $P_5$ into $D$.

\end{lemma}
\begin{proof}
We use the full boundary $S$ throughout. Let $T\subseteq S$ consist of the vertices complete to $D$. A diamond on an edge of $D$ makes $T$ a clique. For a component $E$ of $G-N[D]$, its full exterior boundary is in $S$ and contains some $v\notin T$; otherwise $T$ is a clique cutset separating $E$ from $D$. Choose a transition edge $xy$ in $D$ with $vx=1,vy=0$. A rooted $P_5$ from $v$ into $E$ extends by $y,x$ to a $P_7$, because $D$ misses $E$. Thus $v$ starts no such path. Its matching profile permits \cref{lem:q4-matching-root}, giving $\chi(E)\le4$.

For (1), use four common colors on $D$ and all $E$, and two new colors on $S$. For (2), write $S_g,S_b$ for the good and bad parts on $D$. A component $E$ contacted by $S_b$ has a good contacting root: a root cannot be bad on both anticomplete sets $D,E$, by the two-tail observation. \cref{lem:q4-good-extension} therefore three-colors each such $E$. Give these components and $D$ colors $1,2,3$, and the other components colors $1,2,3,4$.

If $S_b$ is stable, give it $4$ and give $S_g$ colors $5,6$. Otherwise give $S_b$ colors $4,5$ and the stable $S_g$ color $6$. Components using $4$ miss $S_b$, so both palettes are valid across the full boundary.

Apply the first conclusion to the edge $bc$. Its zero graph is exactly $G[A\cup Z]$, proving the four-color bound. The full boundary of a four-chromatic component $D$ of this domain is $B_D\cup C_D$, together with $a$ if $D$ meets $A$; $a$ is isolated in that boundary. Part (1) and criticality force $B_D\cup C_D$ to be nonbipartite.

Every root of this union starts no $P_5$ into $D$, by the apex prefix. Fix the indicated induced four-critical $J\subseteq D$. If a contacting root missed $J$, connectedness would put $J$ in a missed component with a good attachment. \cref{lem:q4-good-extension} would three-color that component and hence $J$, a contradiction. Thus all these roots contact the same $J$.
\end{proof}

\begin{lemma}
\label{cor:c5-zero-shapes}
Every zero component of an induced $C_5$ or $C_7$ is a star or a triangle; singletons and edges count as stars.

\end{lemma}
\begin{proof}
Let $D$ be a zero component that is neither, and let $S$ be its nonempty full boundary. Every member of $S$ contacts the cycle. Apply \cref{lem:bipartite-boundary} to the cycle to obtain $\chi(D)\le4$, and partition $S=S_g\cup S_b$ by goodness on $D$.

A bad root on $D$ is good on the cycle. Its cycle profile is therefore a matching vertex cover containing an edge, by \cref{lem:q4-good-extension}. Two bad roots cannot be adjacent: a neighbor of an edge-completing root misses both ends of that edge, contradicting its own edge coverage. Hence $S_b$ is stable.

If $S_g=\varnothing$, \cref{lem:bipartite-boundary}(1) six-colors $G$. Otherwise \cref{lem:q4-good-extension} gives $\chi(D)\le3$, so it remains to show that $S_g$ is bipartite.

If $D$ contains a triangle, all good roots attach at its same vertex, and $S_g$ is in a matching neighborhood. If $D$ is triangle-free, every good profile covers its edges. A profile containing an edge is isolated in $S_g$, by the completed-edge argument above. Every other profile is a stable vertex cover, hence a bipartition side of connected $D$. No such profile exists when $D$ is nonbipartite. Otherwise both sides have size at least two, since $D$ is not a star, and triangle sparsity makes roots with the same side profile nonadjacent. Thus $S_g$ is bipartite in every case. \cref{lem:bipartite-boundary}(2) again contradicts criticality.
\end{proof}

\subsection{Critical profiles}

\begin{lemma}
\label{lem:critical-profile-toolkit}
\label{lem:missed-triangle}
\label{lem:two-root-contact}
\label{lem:small-profile-critical}
\label{lem:opposite-contact}
\label{lem:profile-edge-no-triangle}
In this lemma, a no-$P_5$ root on $J$ means an exterior vertex starting no rooted $P_5$ into $J$. Let $J$ be induced and $4$-vertex-critical in a $\{P_7,\mathrm{diamond}\}$-free ambient graph of clique number at most three; in particular $\delta(J)\ge3$.

(a) For any contacting no-$P_5$ root $q$, a triangle-containing component of $J-N_J(q)$ is exactly a triangle.

(b) Two contacting no-$P_5$ roots with disjoint profiles each contact every triangle of $J$, have at least three contacts, and contact every component of $J-N_J(q')$ for the other root $q'$.

(c) Let $q,t$ be no-$P_5$ roots, with $t$ contacting every triangle of $J$. If $fg$ is an edge in $N_J(q)$ missed by $t$, neither endpoint lies on a triangle of $J$. No disjointness of the two profiles is required here.

\end{lemma}
\begin{proof}
For (a), the full boundary of a missed component $L$ consists of good roots. If $L$ contains a triangle $T$ and is larger than $T$, \cref{lem:q4-good-extension} gives a common attachment $z\in T$. The other two vertices have no neighbors in $L-T$ or in its boundary. Their degree in $J$ is two, contrary to \cref{lem:critical-core}.

For (b), write the two profiles as $R,R'$. We first prove triangle contact. Suppose $q$ misses $T=\{p,x,y\}$. By (a), $T$ is a component of $J-R$. Root $q'$ cannot miss it too, for its full boundary would then be contained in $R\cap R'=\varnothing$. Say $q'p=1$.

The component $E$ of $J-R'$ containing $xy$ is triangle-free: otherwise (a) gives a second completer of $xy$ besides $p$. Its good boundary profiles cover $xy$. Any boundary vertex outside $T$ would belong to both $R$ and $R'$, because $T$ is a component of $J-R$. Thus this boundary is exactly $\{p\}$. Three-color the proper induced graph $J-E$ by criticality and extend across $E$ using \cref{lem:q4-good-extension}, agreeing at $p$. The contradiction, reflected, proves that both roots contact every triangle. In particular, $J-R$ and $J-R'$ are triangle-free.

If $|R|\le2$, color $R$ distinctly with $1,2$. Each component of $J-R$ has nonempty good boundary in $R$. With one boundary root, extend by \cref{lem:q4-good-extension}, agreeing at that root. With two, apply the \hyperlink{rule:two-boundary-colors}{boundary-coloring observation}: goodness covers every edge, and triangle sparsity makes the common-neighbor set stable. This three-colors $J$, a contradiction. Reflecting gives $|R|,|R'|\ge3$.

Finally suppose $q'$ misses a component $D$ of $J-R$. A shortest path to it is $q',z,w,x$, with $w\in R,x\in D$: disjoint profiles force length at least three, and the no-$P_5$ condition forces length at most three. Vertex $z$ misses $D$. A transition edge in $D$ would give $q',z,w,x',y'$, so $w$ completes $D$. Its matching neighborhood makes $D$ a singleton or an edge. An edge would form a triangle missed by $q'$, already excluded.

For a singleton $x$, its neighborhood is in $R$ and is stable, since $q'$ misses every triangle through $x$. If some $y\in N_J(x)$ misses $z$, then $q',z,w,x,y$ is induced: $q'$ misses $R\cup D$, $zx=0$, and $wy=0$. Hence $N_J(x)\subseteq N_J(z)$, contrary to \cref{lem:critical-core}. Reflect for the other component-contact assertion.

For (c), suppose $fxy$ is a triangle, with $tx=1$. Triangle sparsity gives $qx=qy=gx=gy=0$. By (a) and uniqueness of the completer of $xy$, its component $E$ in $J-N_J(q)$ is triangle-free. Its good boundary profiles cover $xy$, so $g$ misses $E$. If $gz=1$ outside $N_J(q)$, then $z\notin E$ and $zx=zy=zf=0$, the last nonedge following from triangle $qfg$. According as $tz=1$ or $0$, use
\[
t,z,g,f,y\quad\hbox{or}\quad t,x,f,g,z.
\]
Both are rooted $P_5$'s, a contradiction. Thus $N_J(g)=\{f\}$, contrary to $\delta(J)\ge3$. Interchange $f,g$ to finish.
\end{proof}

\subsection{Three colors on each owner domain}

\begin{theorem}
\label{thm:owner-domain-three}
Every owner domain at every triangle of $G$ is three-colorable.

\end{theorem}
\begin{proof}
By \cref{lem:bipartite-boundary}, a four-chromatic domain supplies an induced four-critical $J\subseteq D\subseteq C\cup Z$ and a nonbipartite $A_D\cup B_D$ whose roots contact this same $J$ and start no $P_5$ into it. This union is triangle-free, since an owner edge already has its apex as completer. A shortest odd cycle therefore contains an induced $r,u,v,s$, with $u,v\in B_D$ and $r,s\in A_D$.

Put $R=N_J(u)$, $V=N_J(v)$, $K=J-V$, and $L=J-(R\cup V)$. Since $b$ completes $uv$, the profiles $R,V$ are disjoint. By \cref{lem:critical-profile-toolkit}(b), both roots contact every triangle, have at least three contacts, and contact every component of the other's profile complement. By part (c), no endpoint of a profile edge of $v$ or $r$ lies on a triangle of $J$: the adjacent triangle-contacting root $u$ misses that edge.

\emph{Bipartite complements.} Suppose $K$ contains an odd cycle. A shortest one $F$ has length five or seven, by triangle contact and $P_7$-freeness. Root $r$ is good on $K$ by the prefix $v,b,a$, and contacts $F$, or $r,a,b,v$ is a zero $P_4$, contrary to \cref{cor:c5-zero-shapes}. Its matching vertex cover on $F$ contains exactly one edge $fg$: contact runs have length at most two, missed runs length one. Every $y\in V$ likewise contacts $F$ (or $a,b,v,y$ is a zero $P_4$) and is good there by prefixing $v$. Their completed edges are distinct and avoid $f,g$, whose endpoints lie on no $J$-triangle. Since $|V|\ge3$, $F$ cannot be $C_5$. For $C_7$, label $fg=01$, so $N_F(r)=0135$. Eligible edges are $23,34,45,56$. A root $y$ completing $34$ has profile $1346$, misses $r$, and yields the rooted $P_5$ $r,0,6,y,4$. Reflection $i\mapsto1-i$ excludes $45$, again leaving only two edges. Thus $K$ is bipartite. Swapping $u,v$ and $r,s$ also makes $J-R$ bipartite. Every component of $K$ meets $R$, by \cref{lem:critical-profile-toolkit}(b).

\emph{Fix one palette.} All vertices of $V$ are good on $K$, by prefixing $v$. Endpoints of an edge $pq$ in matching $J[V]$ have stable $J$-neighborhoods by \cref{lem:critical-profile-toolkit}(c); on each nontrivial $K$-component a contacting endpoint therefore sees a whole bipartition side. Their profiles are disjoint. Moreover
\begin{equation}\label{eq:profile-crossing}
x\in R\cap N_K(p),\ z\in L\cap N_K(q)\quad\Longrightarrow\quad xz=1,
\end{equation}
or $u,x,p,q,z$ is a rooted $P_5$. This rule is symmetric in $p,q$.

Call a $K$-edge component wholly in $R$ exceptional. Call a bipartition side pure $R$ or pure $L$ when it is contained in $R$ or $L$, respectively. Fix a two-coloring $\phi$ of $K$ before considering any $V$-edge: on other nontrivial components give a pure $R$-side color $1$ and a pure $L$-side color $2$, orienting freely if neither side is pure. Every component meets $R$; since $R$ is a matching, both sides can be pure $R$ only on an exceptional edge. Thus the prescriptions are consistent. Color singleton components $1$ and exceptional edges arbitrarily.

\emph{Every $V$-edge has a usable endpoint.} We claim that one endpoint has monochromatic $K$-neighborhood under this same $\phi$; an empty neighborhood qualifies. If $px=1$ on an exceptional edge $xy$, triangle $uxy$ gives $py=rx=ry=0$. The sequence $r,a,b,v,p,x,y$ has only the possible chord $rp$, so $rp=1$. Triangle $vpq$ then gives $rq=0$, preventing $q$ from contacting any exceptional edge by the same argument. \eqref{eq:profile-crossing} gives $N_K(q)\subseteq R$, hence its contacted sides are nonexceptional pure $R$-sides, all colored $1$.

Otherwise suppose both neighborhoods are bichromatic. Neither endpoint contacts an exceptional component, and each contacts at least two components. Comparing a contacted side with the mate's contact in another component in \eqref{eq:profile-crossing} makes that side pure $R$ or pure $L$. Thus $p$ sees an $R$-side in $E$ and an $L$-side in $F$. The rule confines $q$ to $E,F$, on their opposite sides; applying it in reverse confines $p$ too. Writing $R_E=E\cap R$, $L_E=E\cap L$, and similarly for $F$, we obtain nontrivial $E,F$ with
\begin{equation}\label{eq:crossed-sides}
N_K(p)=R_E\cup L_F,\quad N_K(q)=L_E\cup R_F.
\end{equation}

Let $H$ be the component of $J-R$ containing $L_E,q,p,L_F$. Its full boundary consists of good roots in $R$, by prefixing $u$. Their profiles cover every edge of bipartite $H$; coverage of $pq$ and \eqref{eq:crossed-sides} make this boundary exactly $R_E\cup R_F$.

These are stable, anticomplete classes, each with identical neighborhoods in $J-H$: all its $K$-neighbors are in $H$, while an outside $V$-neighbor misses $L_E$ or $L_F$ and its good edge-covering profile must contain the whole contacted side. Three-color $J-H$ by criticality, copying a representative's color throughout each class. If their colors agree, two-color $H$ with the other colors. Otherwise call the class colors $1,2$ and use the \hyperlink{rule:two-boundary-colors}{boundary-coloring observation}. Its remaining hypothesis holds: a vertex seeing both classes is in $V$, and such vertices are stable, since adjacent $V$-vertices have disjoint full-side profiles on $E$ and cannot both meet $R_E$. Both alternatives three-color $J$, a contradiction proving the claim.

Finally, for each $V$-edge color a usable endpoint oppositely to its monochromatic neighborhood in $\{1,2\}$ and its mate $3$, choosing either color for an empty neighborhood. Give isolated $V$-vertices color $3$. The selected endpoints are mutually nonadjacent because $J[V]$ is a matching. This extends the fixed $\phi$ to a three-coloring of $J$, the final contradiction.
\end{proof}

\subsection{Owner-cycle boundaries}

We use owner cycles
\[
H=h_0\cdots h_{2k}h_0,\quad k\in\{2,3\},\quad
h_0,h_1,h_3,\ldots,h_{2k-1}\in A,
\]
with all other positions in $Z$. A profile class below contains \emph{all} vertices of $B\cup C$ with that profile; juxtaposition denotes a vertex set.

\begin{lemma}
\label{lem:cycle-profiles}
\label{lem:owner-five-profiles}
\label{lem:seven-profile-palette}
\label{lem:cycle-boundary-coloring}
\label{lem:owner-five-cycle}
\label{lem:owner-seven-cycle}
Contacting profiles on $C_5=r,s,x,t,y,r$ are
\[
rx,rt,st,sy,xy;\quad rxt,rxy,sxy,sty.
\]
On $C_7$, with indices modulo seven, they are
\[
S_i=\{i,i+2,i+4\},\quad E_i=\{i,i+1,i+4\},
\]
or their complements. Stable-profile classes are stable; a profile containing an edge occurs at most once. On $C_5$ the pair roots are bad and the triple roots good. On either cycle, larger-profile roots are mutually nonadjacent, and every contacted neighbor of one has its complementary profile.

If some $B\cup C$ vertex misses $H$, then $\chi(B\cup C)\le3$. For $C_7$ this bound always holds.

\end{lemma}
\begin{proof}
\emph{Profiles and compatible neighbors.} The apex prefix excludes rooted $P_5$'s, hence three consecutive misses; triangle sparsity excludes three consecutive contacts. On $C_5$, the remaining profiles are stable pairs and complementary triples, with $rst$ excluded because $a$ already completes $rs$. On $C_7$ there are three or four contacts: a triple has cyclic gaps $2,2,3$ or $1,3,3$, and a four-contact profile is its complement. Goodness is equivalent to having a stable missed set.

A common contact of adjacent roots forces both neighboring cycle positions to be missed by both roots. Their profiles intersect in at most one position. If their sizes sum to at least $|H|$, an intersection would make their union have size at least $|H|-1$ and at most $|H|-2$. Thus they have disjoint profiles. This proves the complementary-neighbor and larger-profile assertions. Class stability and uniqueness follow from triangle sparsity.

\emph{A missed root yields a palette for the full boundary.} Let $B_0,C_0$ be the parts missing $H$. They cannot both be nonempty: for $e\in B_0,f\in C_0$, the sequence $e,b,c,f$ induces a zero $P_4$ or $C_4$, contrary to \cref{cor:c5-zero-shapes}. Reflecting the owners if necessary, assume $B_0\ne\varnothing$ and $C_0=\varnothing$.

Every bad $C$ root sees all of $B_0$; otherwise its bad tail extends by $v,b,c$ for a missed $B$ root $v$. A good $C$ root $p$ missing such a $v$ has one of the profiles
\[
P_+=\{h_0,h_2\}\cup T_0,\quad
P_-=\{h_1,h_{2k}\}\cup T_0,\quad
T_0=\{h_3,h_5,\ldots,h_{2k-1}\}.
\]
Indeed, its vertex-cover profile misses an end of $h_0h_1$. If it misses $h_0$, it sees $h_1,h_{2k}$. Any other contacted zero position $w$ gives
\[
w,p,h_{2k},h_0,a,b,v.
\]
This is induced, so $p$ misses every other zero position, and edge coverage forces $P_-$. Reflection gives $P_+$.

Put $M_\pm=H\setminus P_\pm$, the stable bad profiles. Fix $v\in B_0$. From every edge of matching $C$, select an endpoint missed by $v$; triangle sparsity guarantees one. These endpoints have profiles $P_\pm$, form a stable set $D$, and have mates of profiles $M_\pm$. Each selected endpoint misses all of $B_0$: its bad mate sees all of $B_0$, and no vertex sees both ends of an owner edge.

Two-color the full matching $B$ with $1,2$, prescribing $M_-:1,M_+:2$; the prescribed classes are disjoint and stable. Give $D\cap P_-$ color $2$, $D\cap P_+$ color $1$, and the stable set $C-D$ color $3$. The complementary-neighbor rule checks contacted $B$ neighbors, while the preceding argument checks the missed roots $B_0$. This proves the bound whenever a root misses $H$.

\emph{The remaining $C_7$ palette.} Suppose all roots contact $H$. Classes $S_i,S_{i+3}$ cannot share an owner: after rotation, same-owner roots $p:024,q:035$ either have an edge with two common neighbors, their apex and $h_0$, or give the induced path
\[
h_1,h_2,p,d,q,h_5,h_6,
\]
where $d$ is that apex. If every $S_i$ occurred, representatives in order $0,3,6,2,5,1,4$ would alternate two owners around an odd cycle. Relabel an empty class as $S_0$ and color
\[
S_2\cup S_4:1,\quad S_6\cup S_1:2,\quad S_3\cup S_5:3.
\]
Each union is stable because its profiles share two positions.

The only possible triple-profile neighbors of $E_i$ are $S_{i+2},E_{i+2},E_{i-2}$. For $E_0=014$, a neighbor misses $0,1$; according as it sees or misses $4$, it is $S_2$ or $E_2,E_5$. Some $E_i$ is absent because $a$ already completes $h_0h_1$. Since each $E_i$ occurs at most once, the $E_i$ induce paths and isolates. Color these greedily, avoiding the preceding color and the color of $S_{i+2}$. Finally each four-contact root avoids the color of its complementary triple class, which is monochromatic; an absent complementary class imposes no restriction. Larger-profile roots are mutually nonadjacent, so this completes the three-coloring.
\end{proof}

\subsection{Local triangle domination}

\begin{lemma}
\label{lem:local-triangle-domination}
Suppose $H=r,s,x,t,y,r$ is an owner $C_5$, every vertex of $B\cup C$ contacts it, $\chi(B\cup C)=4$, and
\[
X=\{p\in B\cup C:N_H(p)=xy\},\quad |X|\ge2.
\]
There is a triangle $T$ such that
\[
N[Q]\cup\{x,y\}\subseteq N[T].
\]

\end{lemma}
\begin{proof}
Call a root matched if it has a mate in its full owner class. Same-owner mates have disjoint profiles. By \cref{lem:cycle-profiles}, the possible pairs of profiles at the ends of an owner edge form the following graph. Its vertices are exact profiles, and an edge records a possible same-owner matching edge:
\[
rx--st--xy--rt--sy--rx,
\]
with leaves $rxy,sxy,rxt,sty$ at $st,rt,sy,rx$, respectively.

Select exactly one endpoint of every owner edge, in both full owner classes. The selected set is stable within each class. If it had no cross-owner edge, give it color $3$, the remaining $B$ vertices color $1$, and the remaining $C$ vertices color $2$. This would contradict $\chi(B\cup C)=4$. Thus every such endpoint selection contains a cross-owner edge.

\emph{The restrictions imposed by a rare edge.} Call $rt$--$xy$ and $st$--$xy$ rare. Suppose $k:rt,l:xy$ are mates in $B$. Every other $B$ root $p$ misses both. Profiles $sy,sty$ give
\[
c,a,s,p,y,l,k,
\]
whereas profiles $xy,sxy$ give
\[
c,a,r,k,l,x,p.
\]
Both are induced paths, excluding those profiles. Consequently every other owner edge in $B$ is $rx$--$st$ or $st$--$rxy$. Reflecting $r\leftrightarrow s,x\leftrightarrow y$ gives the restriction for a rare $st$--$xy$ edge: every other edge is $rt$--$sy$ or $rt$--$sxy$. In particular, each owner has at most one rare edge.

\emph{Both rare types must occur.} Otherwise reflect so that $st$--$xy$ is absent, and select the matched roots of types $rx,rt,rxy,rxt$. The type graph shows that this selects exactly one endpoint of each owner edge. The preceding observation and \cref{lem:cycle-profiles} force a cross edge $uv$ with $u\in B:rx,v\in C:rt$. Their old mates have types
\[
U:st,sy,sty,\quad V:sy,xy,sxy,\quad Uv=Vu=0.
\]
The nonedges follow because an owner edge already has its apex as completer. If $V:xy$, use $s,a,b,u,v,V,y$. If $V:sxy$, the complementary-neighbor rule gives $UV=0$, and use $x,V,v,r,a,b,U$. Thus $V:sy$.

Now let $p:xy$. Exact profiles and ownership give the following restrictions.

If $p\in B$, the sequences $a,\allowbreak c,\allowbreak V,\allowbreak y,\allowbreak p,\allowbreak x,\allowbreak u$ and $a,\allowbreak c,\allowbreak v,\allowbreak u,\allowbreak x,\allowbreak p,\allowbreak y$ have only possible chords $pV$ and $pv$, respectively. If $p\in C$, the sequence $a,\allowbreak b,\allowbreak u,\allowbreak x,\allowbreak p,\allowbreak y,\allowbreak V$ has only possible chord $pu$.

When $p\in B$, $pu=0$ because $u$ has mate $U$. Both chords are therefore forced, making $p$ see $v,V$ and miss $c$ in the triangle $cvV$, impossible. When $p\in C$, $pv=pV=0$ because $vV$ is an owner edge, so $pu=1$. At most one such $p$ can complete $ux$. Hence $|X|\le1$, a contradiction.

\emph{The two rare types determine a replacement triangle.} They occupy opposite owners, since each owner has at most one rare edge. Their restrictions make the matched $rt,st$ roots another endpoint selection. Thus there are $u\in B:rt,v\in C:st$ with $uv=1$. Let $U,V$ be their mates; again $Uv=Vu=0$.

The possibilities for $U$ are $sy,xy,sxy$. Either of the latter two gives $x,U,y,r,a,c,v$, so $U:sy$. Reflection gives $V:rx$. The rare restrictions therefore put $st$--$xy$ in $B$ and $rt$--$xy$ in $C$. Their profile exclusions also show that no $B$ root has type $rx$ and no $C$ root has type $sy$.

Set $T=\{u,v,t\}$. We check its domination of the full old neighborhood, without making any minimum-zero-set assumption in this lemma.

A $B$ root $p$ missing $T$ and seeing $y$ gives
\[
a,c,v,u,U,y,p;
\]
here $pU=0$ because $uU$ is an owner edge. If $p$ misses $y$ as well, its only possible profile is the excluded $rx$. Reflection excludes all missed $C$ roots.

A missed $A$ root $p$ misses $r,s,t$: $r,s$ are mates and $pt=0$. If it sees $x$, use $c,b,u,r,s,x,p$; if it sees $y$, use $b,c,v,s,r,y,p$. Hence it misses all of $H$. The two exhaustive remaining contact cases are
\[
\begin{array}{c|c}
pV=1 & y,t,u,b,c,V,p\\
pV=0 & x,V,v,u,b,a,p.
\end{array}
\]
These are induced: $V:rx$ misses $y,t$, $u:rt$ misses $x,y$, $v:st$ misses $x$, $uV=0$, and $p$ misses $H,u,v,b,c$. The second row also uses $pV=0$.

In all preceding paths, exact profiles remove cycle chords, owner matchings remove same-owner chords, and ownership removes apex chords. The explicitly forced or split contacts are the only remaining ones. All full owner classes are therefore covered by $N[T]$. The old apices are covered because $at=bu=cv=1$, and $x,y$ both see $t$. This proves $N[Q]\cup\{x,y\}\subseteq N[T]$.
\end{proof}

\subsection{Replacing the triangle}

\begin{theorem}
\label{thm:clique-three-bound}
There is no $7$-vertex-critical $\{P_7,\mathrm{diamond}\}$-free graph of clique number three.

\end{theorem}
\begin{proof}
Choose $Q$ among all triangles to minimize the size of its global zero set $Z=V(G)\setminus N[Q]$. \cref{thm:owner-domain-three} three-colors every owner domain.

Every opposite owner pair, say $B\cup C$, has chromatic number four. It is a union of two matchings, hence four-colorable. If three colors $4,5,6$ sufficed, use colors $1,2,3$ on $A\cup Z$ and assign $(a,b,c)=(4,1,2)$. This would six-color $G$.

Partition $Z$ into stable sets $Z_1,Z_2,Z_3$. If $A\cup Z_1,B\cup Z_2,C\cup Z_3$ were all bipartite, give them disjoint palettes $\{1,2\},\{3,4\},\{5,6\}$, respectively, and assign $(a,b,c)=(3,5,1)$. This again six-colors $G$. After permuting owners, $A\cup Z_1$ therefore contains a shortest odd cycle. It is triangle-free: $Z_1$ is stable, and an owner edge already has its apex as completer. Its length is five or seven. Zero positions are isolated and owner runs have length at most two, so it is an owner cycle of the form used above.

\cref{lem:cycle-profiles} excludes $C_7$, since $\chi(B\cup C)=4$. The same lemma forces every $B\cup C$ root to contact the remaining $C_5=r,s,x,t,y,r$.

To apply the local domination lemma, we first obtain its hypothesis $|X|\ge2$ from the globally chosen $Q$. The triangle $\{a,r,s\}$ covers every old apex and owner vertex except possibly the profile class $X$ of type $xy$, and also covers $x,y$. Its global zero set is therefore contained in
\[
(Z\setminus\{x,y\})\cup X.
\]
Minimality of $|Z|$ among all triangles gives $|Z|\le |Z|-2+|X|$, so $|X|\ge2$.

\cref{lem:local-triangle-domination} now supplies a triangle $T$ with
\[
N[Q]\cup\{x,y\}\subseteq N[T].
\]
Taking complements in the same ambient $V(G)$, its global zero set satisfies
\[
Z_T\subseteq Z\setminus\{x,y\},\quad
|Z_T|\le |Z|-2<|Z|.
\]
This contradicts the choice of $Q$.
\end{proof}

\section{\texorpdfstring{Proof of \cref{thm:main}}{Proof of the main theorem}}\label{sec:main-proof}
\begin{proof}[Proof of \cref{thm:main}]

The null and nonnull edgeless graphs satisfy the stated conventions. If the bound failed, \cref{cor:seven-critical-reduction} would give a seven-vertex-critical $\{P_7,\mathrm{diamond}\}$-free graph of clique number three. \cref{thm:clique-three-bound} excludes it. This proves the theorem.
\end{proof}
\section*{Declaration of competing interest}
The authors declare that they have no known competing financial interests or personal relationships that could have appeared to influence the work reported in this paper.

\section*{Acknowledgements}
The research was partially supported by RGC Competitive Earmarked Research (Grant No. 16308821), by the Foundation for Cultivated Young Talents of Fujian Province, China (Grant No. 2026350294), by the Natural Science Foundation of Fujian Province, China (Grant No. 2026J001968), by the Youth Foundation of Fujian Province (Grant No. JZ240035), by the Minnan Normal University Foundation (Grant No. KJ2023002), and by the Fujian Key Laboratory of Granular Computing and Applications (Minnan Normal University), the Institute of Meteorological Big Data-Digital Fujian, and the Fujian Key Laboratory of Data Science and Statistics..

\section*{Data availability}
The proof of the upper bound is contained in the article and does not rely on computational data or computer-generated certificates.

\section*{Declaration of generative AI use}
ChatGPT assisted with the interface analysis in \cref{local-common-missing-label} and the rooted-coloring argument in the proof of \cref{thm:owner-domain-three}. It also assisted with language editing. The authors take full responsibility for the accuracy and content of this work.

\end{document}